\documentclass[11pt]{article}
\usepackage[a4paper,margin=25mm]{geometry}
\usepackage[T1]{fontenc}
\usepackage{lmodern,amsmath,amssymb,amsthm,mathrsfs,microtype,graphicx,tikz,booktabs,longtable,array,enumitem}
\usepackage[colorlinks=true,linkcolor=blue!45!black,citecolor=blue!45!black,urlcolor=blue!45!black]{hyperref}
\usetikzlibrary{arrows.meta,decorations.pathreplacing,calc}
\usepackage[font=footnotesize]{caption}
\numberwithin{equation}{section}
\newtheorem{theorem}{Theorem}[section]
\newtheorem{lemma}[theorem]{Lemma}
\newtheorem{proposition}[theorem]{Proposition}
\newtheorem{corollary}[theorem]{Corollary}
\newtheorem{maintheorem}{Theorem}

\theoremstyle{definition}
\newtheorem{definition}[theorem]{Definition}
\theoremstyle{remark}

\newtheorem{example}[theorem]{Example}
\newcommand{\R}{\mathbb R}
\newcommand{\PP}{\mathbb P}
\newcommand{\EE}{\mathbb E}
\newcommand{\cL}{\mathcal L}
\newcommand{\Rc}{\mathscr R}
\setlist[enumerate]{itemsep=2pt,topsep=3pt}

\title{Accumulating diffusion zeros and long-time dynamics\\
of one-dimensional stochastic differential equations}
\author{Tuan Anh Phan\thanks{Department of Mathematics, Augusta University,
Augusta, Georgia, USA.
Email: \href{mailto:tuphan@augusta.edu}{\texttt{tuphan@augusta.edu}}.}
\and Jinqiao Duan\thanks{Department of Mathematics and Department of Physics,
School of Sciences, Great Bay University, Dongguan, Guangdong 523000, China.
Email: \href{mailto:duan@gbu.edu.cn}{\texttt{duan@gbu.edu.cn}}.}
\and Jianjun Paul Tian\thanks{Department of Mathematical Sciences,
New Mexico State University, Las Cruces, New Mexico, USA. Corresponding author.
Email: \href{mailto:jtian@nmsu.edu}{\texttt{jtian@nmsu.edu}}.}}
\date{}

\begin{document}
\maketitle

\begin{abstract}
We study scalar It\^o equations $dX_t=f(X_t)\,dt+\sigma(X_t)\,dW_t$ with locally
Lipschitz coefficients whose diffusion coefficient may vanish on an arbitrary
closed set. A common zero of $f$ and $\sigma$ is an equilibrium: it is never
reached and it separates the line. A zero of $\sigma$ alone is a one-way gate:
it is crossed at most once, in the direction of the drift. A path may cross
infinitely many accumulating gates in finite time. Between consecutive
equilibria, or out to $\pm\infty$, almost every path either converges to an
outer endpoint or eventually stays in a single recurrent component. We
characterize all stationary probabilities. When there is an interior diffusion
zero, we also determine the weak limit of the law and the almost-sure limit
of the empirical measures. We bound passage times across blocks independently
of the number of zeros and give an exact mean-passage criterion for the stage
between a repelling equilibrium and its first outward gate.
\end{abstract}



\section{Introduction and main results}\label{sec:intro}

\subsection{Two kinds of diffusion zeros}

Consider the scalar It\^o equation
\begin{equation}\label{eq:sde}
 dX_t=f(X_t)\,dt+\sigma(X_t)\,dW_t,\qquad X_0=x,
\end{equation}
with locally Lipschitz coefficients. Where $\sigma\ne0$, the noise moves the
path in both directions and every nearby level is reached. Where $\sigma=0$,
the noise switches off and only the drift acts. Such zeros come in two kinds.

At a common zero $e$ of $f$ and $\sigma$, the process stands still. Such an
\emph{equilibrium} is never reached in finite time from any other state, so
it splits the state space into invariant pieces. A path can only approach it as
$t\to\infty$.

At a zero $h$ of $\sigma$ with $f(h)\ne0$, the drift pushes the path across. The
crossing always goes in the direction of $f(h)$, and it happens at most once.
We call such a point a \emph{gate}. The classical name for this kind of point
in general one-dimensional diffusions is a shunt point~\cite{Li2021}.

When there are finitely many gates, long-time behavior follows by counting.
Each gate is crossed at most once, so after finitely many crossings the path
stays in one interval, where the classical scale and speed theory applies. When
gates \emph{accumulate}, this argument breaks down. The path can cross
infinitely many gates in a finite time, and an accumulation point of gates
is not isolated, so the local classification of isolated singular
points~\cite{ChernyEngelbert} does not reach it.

\begin{figure}[htbp]
\centering
\resizebox{\linewidth}{!}{%
\begin{tikzpicture}[>=Stealth,font=\small]
 \fill[teal!6] (0,-0.32) rectangle (2.6,0.32);
 \fill[teal!6] (12.2,-0.32) rectangle (14.4,0.32);
 \fill[teal!20] (8.0,-0.32) rectangle (10.4,0.32);
 \draw[line width=0.9pt] (0,0) -- (14.4,0);
 \draw[line width=0.9pt,->] (14.4,0) -- (15.0,0);
 \node[right] at (15.0,0) {$r=+\infty$};
 \fill (0,0) circle (2.6pt);
 \node[below=4pt] at (0,-0.1) {$l$};
 \node[above=4pt,align=center,font=\scriptsize] at (0,0.25) {equilibrium};
 \foreach \x/\s in {2.6/-1,4.6/1,10.4/-1,12.2/1}{
   \draw[line width=0.8pt] (\x,-0.2) -- (\x,0.2);
   \draw[->,line width=0.8pt] (\x,0.45) -- ++(0.42*\s,0);
 }
 \foreach \x in {6.0,6.7,7.15,7.45,7.65,7.78,7.86,7.92,7.96,7.985}{
   \draw[line width=0.5pt] (\x,-0.13) -- (\x,0.13);
 }
 \draw[line width=0.8pt] (8.0,-0.2) -- (8.0,0.2);
 \node[below=3pt] at (8.0,-0.18) {$c$};
 \draw[->,line width=0.8pt] (6.0,0.45) -- (8.0,0.45);
 \node[above,font=\scriptsize] at (7.0,0.48) {$f\ge\beta>0$};
 \draw[decorate,decoration={brace,mirror,amplitude=4pt}] (6.0,-1.15) -- (8.0,-1.15)
   node[midway,below=5pt,align=center,font=\scriptsize]
   {accumulation block:\\ crossed in finite mean time};
 \node[below=8pt,align=center,font=\scriptsize] at (1.3,-0.2)
   {converge to $l$,\\ or recurrent};
 \node[below=8pt,align=center,font=\scriptsize] at (3.6,-0.2)
   {two exits:\\ branch by scale};
 \node[below=8pt,align=center,font=\scriptsize] at (5.3,-0.2)
   {pass\\ right};
 \node[below=8pt,align=center,font=\scriptsize] at (9.2,-0.2)
   {positive\\ recurrent};
 \node[below=8pt,align=center,font=\scriptsize] at (11.3,-0.2)
   {two exits:\\ branch by scale};
 \node[below=8pt,align=center,font=\scriptsize] at (13.5,-0.2)
   {escape to $+\infty$,\\ or recurrent};
 \node[above=12pt,font=\scriptsize] at (2.6,0.2) {gate};
 \node[above=12pt,font=\scriptsize] at (4.6,0.2) {gate};
 \node[above=12pt,font=\scriptsize] at (10.4,0.2) {gate};
 \node[above=12pt,font=\scriptsize] at (12.2,0.2) {gate};
\end{tikzpicture}}
\caption{The picture behind Theorem~\ref{thm:A}. The filled dot is an
equilibrium; vertical ticks are gates (zeros of $\sigma$ with $f\ne0$), and the
arrow above a gate shows the only direction in which it can be crossed. Gates
accumulate at $c$ from the left, and the whole block is crossed in finite mean
time. The dark interval is bounded by two inward gates, so it is positive
recurrent. The light intervals are recurrent if their outer endpoint ($l$ or
$+\infty$) is repelling; otherwise paths that enter them converge to that
endpoint. Table~\ref{tab:local} decides which. From any initial state, the path
makes its one-way crossings and then either converges to an outer endpoint or
stays forever in one recurrent interval.}
\label{fig:big-picture}
\end{figure}
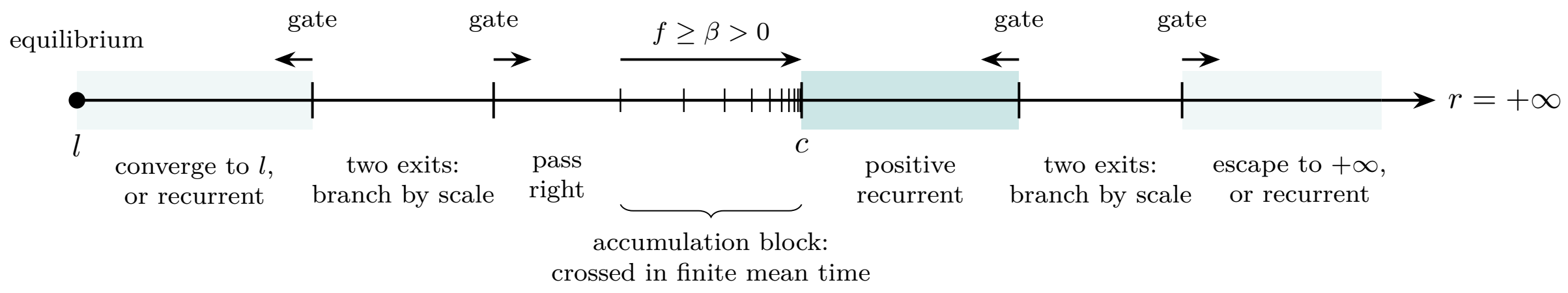

This paper shows that counting is never needed. Between two consecutive
equilibria, or out to $\pm\infty$, almost every path either converges to an
outer endpoint or eventually stays in a single recurrent component, whatever
the zero set of $\sigma$ (Theorem~\ref{thm:A}). This gives all stationary
distributions and, when there is an interior gate, the limiting law
(Theorem~\ref{thm:B}). Theorem~\ref{thm:C} bounds the time to cross a block
independently of the number of gates it contains.
Figure~\ref{fig:big-picture} shows the whole picture.

\subsection{Main results}

Let $-\infty\le l<r\le\infty$. The coefficients are locally Lipschitz on
$(l,r)$ and on a neighbourhood of each finite endpoint. Each finite endpoint is
an equilibrium, and $(f,\sigma)\ne(0,0)$ on $(l,r)$. Let
\[
 Z=\{x\in(l,r):\sigma(x)=0\},
\]
an arbitrary relatively closed subset of $(l,r)$. Every point of $Z$ is a gate.
The \emph{components} are the connected components of $(l,r)\setminus Z$. They
are at most countably many open intervals $J=(u,v)$, and on each of them the
scale function $S_J$ and speed density $m_J$ are defined by~\eqref{eq:scale}
below. An end of $J$ is \emph{attracting} if $S_J$ has a finite limit there and
\emph{repelling} otherwise. Let $\Rc$ be the set of components whose two ends are
both repelling, and $\Rc_+\subset\Rc$ those with finite speed mass
$M_J=\int_Jm_J<\infty$. Solutions are maximal, with lifetime $\zeta$.

\begin{maintheorem}[Global alternative]\label{thm:A}
For every $x\in(l,r)$, the solution stays in $(l,r)$ for $t<\zeta$, and almost
surely exactly one of the following holds.
\begin{enumerate}[label=\textup{(\roman*)}]
\item $X_t\to l$ or $X_t\to r$ as $t\uparrow\zeta$. Convergence to a finite
equilibrium takes infinite time. Explosion, $\zeta<\infty$, can occur only
toward an infinite endpoint.
\item There are a component $J=(u,v)\in\Rc$ and an almost surely finite random
time after which $X_t\in J$ forever. Then $\zeta=\infty$,
$\liminf_{t\to\infty}X_t=u$ and $\limsup_{t\to\infty}X_t=v$.
\end{enumerate}
Every $J\in\Rc$ is invariant, and each of its ends that lies in $Z$ is a gate
whose drift points into $J$. The process on $J$ is positive recurrent if
$J\in\Rc_+$ and null recurrent otherwise.
\end{maintheorem}

Theorem~\ref{thm:A} contains no finiteness or discreteness assumption on $Z$.
Its proof (Section~\ref{sec:accumulation}) excludes an interior limit by local
exit estimates. A path with distinct lower and upper limits must revisit the
levels between them, so no gate can lie strictly between those limits.
When $Z$ is finite, this gives an explicit route through finitely many
intervals, with at most one random branch
(Figure~\ref{fig:route}, Section~\ref{sec:finite}).

Theorem~\ref{thm:B}, proved in Section~\ref{sec:laws}, turns this pathwise
picture into statements about laws. Write $q_l(x)$ and $q_r(x)$ for the
probabilities of the two events in (i). For $J\in\Rc$, write $w_J(x)$ for the
probability of eventual entry into $J$. Let $\pi_J=m_J\,dx/M_J$ for
$J\in\Rc_+$. We write $\mathcal L_x(X_t)$ for the law of $X_t$ under $\PP_x$.
After an explosion, $X_t$ is set equal to the infinite endpoint reached, so
$\mathcal L_x(X_t)$ is a probability measure on the two-point compactification
of $[l,r]$, which adds $\pm\infty$ as points when $l$ or $r$ is infinite. The
arrow $\Longrightarrow$ denotes weak convergence of probability measures on this
compact space: $\int g\,d\mathcal L_x(X_t)\to\int g\,d\mu$ for every continuous
function $g$ on it. Convergence in total variation means that
$\sup_A|\PP_x(X_t\in A)-\mu(A)|\to0$ over all Borel sets $A$.

\begin{maintheorem}[Stationary and limiting laws]\label{thm:B}
\begin{enumerate}[label=\textup{(\roman*)}]
\item The stationary probability measures on the union of $(l,r)$ and the finite
endpoints are exactly the mixtures
\[
 \alpha\,\delta_l+\beta\,\delta_r+\sum_{J\in\Rc_+}\gamma_J\,\pi_J,
\]
where $\alpha,\beta,\gamma_J\ge0$ sum to one, and $\alpha=0$ $(\beta=0)$ when
$l$ $(r)$ is infinite.
\item Suppose $Z\ne\varnothing$. Each null recurrent component touches exactly
one outer endpoint $e_J\in\{l,r\}$, and its speed mass is infinite only there.
For every $x\in(l,r)$, as $t\to\infty$,
\begin{equation}\label{eq:weak-limit}
 \mathcal L_x(X_t)\Longrightarrow
 Q_l(x)\,\delta_l+Q_r(x)\,\delta_r+\sum_{J\in\Rc_+}w_J(x)\,\pi_J,
\end{equation}
where $Q_e(x)=q_e(x)+\sum_{J\in\Rc\setminus\Rc_+,\,e_J=e}w_J(x)$ for
$e\in\{l,r\}$. The convergence holds in total variation if and only if
$Q_l(x)+Q_r(x)=\PP_x(\zeta<\infty)$; without explosion, exactly when
$Q_l(x)=Q_r(x)=0$.
\item Suppose $Z\ne\varnothing$ and $x\in(l,r)$. Almost surely, the empirical
measures $t^{-1}\int_0^t\delta_{X_s}\,ds$ converge weakly to $\delta_e$ on
convergence to the outer endpoint $e$, to $\pi_J$ on eventual entry into
$J\in\Rc_+$, and to $\delta_{e_J}$ on eventual entry into
$J\in\Rc\setminus\Rc_+$.
\end{enumerate}
\end{maintheorem}

A null recurrent component with a unique infinite-speed end keeps visiting
every interior level, yet its law converges to a point mass at that end.
If the end is infinite, mass escapes every compact subset of the real line.
When $Z=\varnothing$ and a recurrent component has infinite speed mass at
both ends, Theorem~\ref{thm:null} locates the support of subsequential limits,
and the weights may oscillate without limit~\cite{Rudnicki1993}.

On the recurrent alternative of Theorem~\ref{thm:A}, entry into the selected
component occurs in finite time. A path converging to an outer endpoint may
cross infinitely many gates without a last crossing. Theorem~\ref{thm:C}
estimates block passage and, for finite routes, the time of the last crossing.
For a level $a$, let $T_a=\inf\{t\in[0,\zeta):X_t=a\}$ be the
first time the path reaches $a$, with $\inf\varnothing=\infty$. When $Z$ is
finite, the path crosses finitely many gates, each at most once
(Section~\ref{sec:finite}). The \emph{selection time} $\tau_*$ is the time of
the last gate crossing, with $\tau_*=0$ if no gate is crossed. After $\tau_*$ the
path stays in the closure of the component in which its route ends.

\begin{maintheorem}[Passage and selection times]\label{thm:C}
\begin{enumerate}[label=\textup{(\roman*)}]
\item \emph{Blocks.} Let $[h,d]\subset(l,r)$, $\sigma(h)=0$ and $f\ge\beta>0$ on
$[h,d]$. For $x\in[h,d)$, $T_d<\zeta$ almost surely,
\[
 \EE_xT_d\le\frac{d-x}{\beta},\qquad
 \sup_{x\in[h,d)}\PP_x\Big(T_d>\frac{2k(d-h)}{\beta}\Big)\le2^{-k}
 \quad(k=0,1,2,\ldots),
\]
whatever the zeros of $\sigma$ in $(h,d)$.
\item \emph{Routes.} Suppose $Z$ is finite and every component visited before
the final one is bounded by two gates. Choose $C>0$ bounding the mean time
to the next exit from each such component, uniformly over its interior starts
and any inward entrance gate.
Then the selection time $\tau_*$ satisfies, for $0<\theta<\log2/(2C)$,
\[
 \EE_xe^{\theta\tau_*}\le\Big(\frac{e^{2C\theta}}{2-e^{2C\theta}}\Big)^{|Z|}.
\]
\item \emph{The stage next to an equilibrium.} Let $l$ be a finite equilibrium
that is repelling from the right. Suppose $Z$ has a smallest point $h$ and
$f(h)>0$. For $l<x<h$, with $S,m$ on $(l,h)$ and $B=S(h-)$,
\[
 \EE_xT_h=(B-S(x))\int_l^xm(y)\,dy+\int_x^h(B-S(y))\,m(y)\,dy .
\]
In particular $\EE_xT_h<\infty$ if and only if the speed mass near $l$ is finite.
\item \emph{Strict repulsion.} If moreover $f,\sigma$ are $C^1$ near $l$ and
$\lambda=f'(l)-\tfrac12\sigma'(l)^2>0$, then for $0<c_0<\lambda$ and $d>l$ close
to $l$,
$\EE_xT_d\le c_0^{-1}\log\big((d-l)/(x-l)\big)$ for $l<x<d$.
\end{enumerate}
\end{maintheorem}

Part~(i) applies even to a block with infinitely many gates. In part~(iii)
the outward gate is hit almost surely. If the speed mass near the repelling
equilibrium is infinite, this first passage already has infinite mean
(Example~\ref{ex:infinite-mean}).

\subsection{Local tests}

The theorems are stated in terms of attracting ends and finite speed mass. Both
can be read off the coefficients. Let $e$ be a finite end of a component $J$
(an equilibrium or a gate). Put $s=+1$ if $J$ lies to the right of $e$ and $s=-1$
otherwise, and measure the distance into $J$ by $y=s(x-e)>0$. In this
\emph{inward coordinate} the drift and the diffusion coefficient become
\[
 \mu(y)=s\,f(e+sy),\qquad \nu(y)=\sigma(e+sy),
\]
so that $\mu>0$ means the drift points into $J$, away from $e$. We expand
\[
 \mu(y)=A\,y^{p}+O(y^{p+1}),\qquad \nu(y)=B\,y^{n}+O(y^{n+1}),\qquad
 k=\frac{2A}{B^2},
\]
with integers $p\ge0$, $n\ge1$ and $AB\ne0$. At an equilibrium $p\ge1$; at a
gate $p=0$. At $+\infty$ we use $x$ itself and write
$f(x)\approx A\,x^{p}$ and $\sigma(x)\approx B\,x^{n}$ as $x\to\infty$, with
real exponents, so that $A>0$ means the drift points toward $+\infty$. The
precise hypotheses are~\eqref{eq:orders} and~\eqref{eq:growth} in
Section~\ref{sec:local}.

\begin{proposition}[Local tests]\label{prop:local}
The ends of a component are classified as in Table~\ref{tab:local}. The rows for
gates require no finite-order expansion, only local Lipschitz continuity.
The proposition is proved at the end of Section~\ref{sec:component}.
\end{proposition}

\begin{table}[htbp]
\centering
\small
\renewcommand{\arraystretch}{1.18}
\begin{tabular}{@{}p{50mm}p{45mm}p{48mm}@{}}
\toprule
End of the component & Attracting (finite scale) & Finite speed mass nearby\\
\midrule
Equilibrium, $C^2$, $\sigma'(e)\neq0$
 & iff $\lambda=f'(e)-\tfrac12\sigma'(e)^2<0$ & iff $\lambda>0$\\
Equilibrium, $p<2n-1$ & iff $A<0$ & iff $A>0$\\
Equilibrium, $p=2n-1$ & iff $k<1$ & iff $k>2n-1$\\
Equilibrium, $p>2n-1$ or $\mu\equiv0$ & always & never\\
Gate, drift into the component & never & always\\
Gate, drift out of the component & always; hit in finite time on approach & ---\\
$+\infty$, $p>2n-1$ & iff $A>0$ & iff $A<0$\\
$+\infty$, $p=2n-1$ & iff $k>1$ & iff $k<2n-1$\\
$+\infty$, $p<2n-1$ & never & iff $n>\tfrac12$\\
\bottomrule
\end{tabular}
\caption{Local tests (Proposition~\ref{prop:local}). The first row also covers
$f'(e)=0$ when $\sigma'(e)\ne0$. The end $-\infty$ follows by the reflection $x\mapsto-x$. Under $C^1$
regularity alone, $\lambda=0$ is undecided (Example~\ref{ex:critical}). A flat
zero, such as $\sigma(x)=e^{-1/x^2}$, is not decided by any finite order.}
\label{tab:local}
\end{table}

At an equilibrium, attracting from one side means stable in probability from
that side (Lemma~\ref{lem:scale}). For instance, the classical threshold
$\lambda<0$ for $dX=aX\,dt+bX\,dW$ reads $a<b^2/2$. The two columns answer
different questions: for $n\ge2$ and $1\le k\le2n-1$, an equilibrium is
repelling, yet its speed mass is infinite. This is the mechanism behind null
recurrence in Theorem~\ref{thm:B}.

\subsection{Related work}

The analytic tools are classical. Scale functions and speed measures classify
the boundary behavior of a regular diffusion on one
interval~\cite{Feller1952,Mandl1968,Karlin}, and Feller's test decides
explosion~\cite{KaratzasShreve}. Stability in probability and its Lyapunov
criteria are developed in~\cite{Khasminskii2012}. A solution with Lipschitz
coefficients never reaches the common zero set of all its coefficients; this is
proved in any dimension in~\cite{Lyons}. For degenerate scalar equations,
Engelbert and Schmidt characterized existence and uniqueness in terms of the
zeros of $\sigma$~\cite{EngelbertSchmidt}. Cherny and Engelbert gave a complete
local classification of isolated singular points~\cite{ChernyEngelbert}. One-way
points of general one-dimensional diffusions are the shunt points
of~\cite{Li2021}. For a single recurrent component, the empirical-measure
conclusion when precisely one end has infinite speed mass also appears in
\cite[Theorem~2.11]{Lemaire2007}. Here we also identify the limiting
fixed-time law and combine the component laws using their entry probabilities.
When both ends have infinite speed mass, the law need not converge at all:
Rudnicki~\cite{Rudnicki1993} gives sufficient conditions for convergence of the
weights at $\pm\infty$ and an example in which they oscillate forever.

Our setting differs in two ways. First, the zeros of $\sigma$ need not be
isolated. They may accumulate at gates, where local classification does not
apply and counting fails. Second, we work globally, between consecutive
equilibria or out to infinity. There we obtain a pathwise alternative
(Theorem~\ref{thm:A}), all stationary laws and the weak limit when an interior
gate is present (Theorem~\ref{thm:B}), and bounds on block passage that do not
depend on the number of gates crossed (Theorem~\ref{thm:C}).

\subsection{Organization}

Section~\ref{sec:setting} fixes the setting and terminology and collects them in
a glossary. Section~\ref{sec:local} proves the local tests of
Table~\ref{tab:local} and the behavior of a single component.
Section~\ref{sec:finite} treats finitely many gates.
Section~\ref{sec:accumulation} proves Theorem~\ref{thm:A}.
Sections~\ref{sec:laws} and~\ref{sec:times} prove Theorems~\ref{thm:B}
and~\ref{thm:C}. Section~\ref{sec:discussion} lists open problems.
Classical tools are collected in the appendix.

\section{Setting, conventions and glossary}\label{sec:setting}

This section fixes the standing assumptions, the two kinds of diffusion zeros,
and the scale and speed functions on which all later arguments rest. Its only
new statement is the elementary fact that equilibria are never reached
(Proposition~\ref{prop:inaccessible}). Table~\ref{tab:glossary} collects the
terminology.

\subsection{Standing assumptions}

Throughout, $-\infty\le l<r\le\infty$ and
\begin{equation}\label{eq:standing}
\begin{gathered}
 f,\sigma\ \text{are locally Lipschitz on } (l,r)
 \text{ and on a neighbourhood of each finite endpoint},\\
 f(e)=\sigma(e)=0\ \text{at each finite endpoint } e,\qquad
 (f(x),\sigma(x))\neq(0,0)\quad(l<x<r).
\end{gathered}
\end{equation}
Thus $l$ and $r$ are consecutive equilibria, or $\pm\infty$. For each
$x\in(l,r)$, equation~\eqref{eq:sde} has a unique maximal strong solution with
lifetime $\zeta\in(0,\infty]$, defined on a filtered probability space with a
Brownian motion $W$~\cite{Ikeda}. We write $\PP_x$ and $\EE_x$ for the
corresponding law and expectation. For a level $a$ and an open set $D$,
\[
 T_a=\inf\{t\in[0,\zeta):X_t=a\},\qquad
 \tau_D=\zeta\wedge\inf\{t\in[0,\zeta):X_t\notin D\},
\]
with $\inf\varnothing=\infty$. The generator is
$\cL V=fV'+\tfrac12\sigma^2V''$.

\paragraph{Conventions and scope.}
All results concern the maximal solution, so no global non-explosion assumption
is made. ``Almost surely'' always refers to $\PP_x$ for a fixed initial state
$x$. Stability and attraction are one-sided notions, taken from the side of the
component under study. Convergence to a finite equilibrium always means
convergence as $t\to\infty$, never hitting. After an explosion, $X_t$ is set
equal to the infinite endpoint reached, for $t\ge\zeta$. The local results of
Section~\ref{sec:local} hold on any interval on which their hypotheses hold;
the global results need~\eqref{eq:standing} only.

\subsection{Equilibria, gates and components}

\begin{definition}\label{def:points}
Let $x\in(l,r)$ or let $x$ be a finite endpoint.
\begin{enumerate}[label=\textup{(\roman*)}]
\item $x$ is an \emph{equilibrium} if $f(x)=\sigma(x)=0$.
\item $x$ is a \emph{gate} if $\sigma(x)=0$ and $f(x)\neq0$. It is
\emph{rightward} if $f(x)>0$ and \emph{leftward} if $f(x)<0$.
\item The \emph{components} are the connected components of
$(l,r)\setminus Z$, where $Z=\{x\in(l,r):\sigma(x)=0\}$.
\end{enumerate}
\end{definition}
Under~\eqref{eq:standing} there are no equilibria inside $(l,r)$, so every
point of $Z$ is a gate. In particular every accumulation point of $Z$ inside
$(l,r)$ is a gate; by continuity, its drift has the same sign on a whole
neighbourhood. The set $Z$ may be infinite, non-discrete, or uncountable.

\begin{proposition}[Equilibria are never reached]\label{prop:inaccessible}
Let $e$ be an equilibrium and $x\neq e$. Then
$\PP_x\big((X_t-e)(x-e)>0\ \text{for all } t<\zeta\big)=1$.
\end{proposition}

\begin{proof}
Work with the given extension across $e$ if $e$ is a finite endpoint.
Since $f(e)=\sigma(e)=0$, local Lipschitz continuity implies that
$a(z):=\frac{f(z)}{z-e}$ and $b(z):=\frac{\sigma(z)}{z-e}$
are bounded on compact sets, with arbitrary values assigned at $e$.
Before the first hitting time $T_e$, It\^o's formula gives
\[
 \log|X_t-e|
 =\log|x-e|
 +\int_0^t\!\left(a(X_s)-\tfrac12 b(X_s)^2\right)ds
 +\int_0^t b(X_s)\,dW_s .
\]
If $T_e$ were finite, the left-hand side would tend to $-\infty$ as
$t\uparrow T_e$. The right-hand side has a finite limit almost surely:
after localization to a compact interval, the drift integrand is bounded
and the stochastic integral extends continuously to $T_e$.
Thus $X$ never hits $e$, and continuity keeps $X_t-e$ of the same sign
as $x-e$ for every $t<\zeta$. The multidimensional inaccessibility
result is proved in~\cite{Lyons}.
\end{proof}

The companion statement for gates, that a gate is crossed at most once and only
in the direction of its drift, is Theorem~\ref{thm:one-way} in
Section~\ref{sec:local}.

\subsection{Scale, speed and stability}

On a component $J=(u,v)$, fix $c\in J$ and set
\begin{equation}\label{eq:scale}
 \rho_J(x)=\exp\Big(-\int_c^x\frac{2f(z)}{\sigma^2(z)}\,dz\Big),\qquad
 S_J(x)=\int_c^x\rho_J(y)\,dy,\qquad
 m_J(x)=\frac{2}{\sigma^2(x)\rho_J(x)}.
\end{equation}
$S_J$ is the \emph{scale function} and $m_J$ the \emph{speed density}. A
different $c$ changes $S_J$ by an increasing affine map and $m_J$ by the
reciprocal factor, so finiteness of the limits $S_J(u+)$, $S_J(v-)$ and of
the speed mass $M_J=\int_Jm_J$ does not depend on $c$. Scale is never continued
across a zero of $\sigma$. For $u<a<x<b<v$, the scale function gives the
two-barrier formula, which is the probability of reaching $b$ before $a$,
\begin{equation}\label{eq:two-barrier}
 \PP_x(T_b<T_a)=\frac{S_J(x)-S_J(a)}{S_J(b)-S_J(a)}.
\end{equation}

\begin{definition}[Attracting and repelling ends]\label{def:ends}
The left end $u$ of $J$ is \emph{attracting} if $S_J(u+)>-\infty$ and
\emph{repelling} if $S_J(u+)=-\infty$. The right end is attracting if
$S_J(v-)<\infty$ and repelling otherwise.
\end{definition}

Letting $a\downarrow u$ in~\eqref{eq:two-barrier} shows the meaning: at a
repelling end, every level inside $J$ is reached before the path gets close to
that end.

\begin{definition}[Stability in probability]\label{def:stable}
An equilibrium $e$ is \emph{stable in probability from the right} if for every
$\varepsilon>0$ and $\eta\in(0,1)$ there is $\delta>0$ such that
\[
 e<x<e+\delta\quad\Longrightarrow\quad
 \PP_x\big(\sup_{t\ge0}|X_t-e|\ge\varepsilon\big)<\eta .
\]
It is \emph{asymptotically stable in probability from the right} if, in
addition, $\PP_x(X_t\to e)\to1$ as $x\downarrow e$. Stability from the left is
defined symmetrically, and $e$ is \emph{unstable} from a side if it is not
stable from that side.
\end{definition}

For an equilibrium end with $\sigma\neq0$ on the adjacent side, these notions
coincide with attraction: $e$ is stable from that side if and only if it is
asymptotically stable, if and only if it is an attracting end
(Lemma~\ref{lem:scale}). The scale function also governs recurrence; the
following lemma is classical~\cite[Ch.~15]{Karlin}, \cite{MeynTweedie1993}.

\begin{lemma}[Recurrence and speed]\label{lem:recurrence}
Suppose both ends of a component $J$ are repelling and the process started in
$J$ stays in $J$. Then it is recurrent: every level of $J$ is hit, and any two
levels are visited alternately infinitely often. It is positive recurrent if
and only if $M_J<\infty$. Then $\pi_J(dx)=m_J(x)\,dx/M_J$ is its unique
stationary probability on $J$; for each $x\in J$, the law of $X_t$ converges to
$\pi_J$ in total variation, and time averages of bounded Borel functions
converge almost surely to their $\pi_J$-means. If $M_J=\infty$, it is null
recurrent and has no stationary probability on $J$.
\end{lemma}

The density $\pi_J$ is the zero-flux solution of the stationary Fokker--Planck
equation, $f\pi_J-\tfrac12(\sigma^2\pi_J)'=0$. A formal stationary solution with
infinite mass does not define a stationary probability.

\subsection{Glossary}

Table~\ref{tab:glossary} explains the terms used in this paper in plain
language. Precise definitions are given where each term first appears.

{\small
\renewcommand{\arraystretch}{1.14}
\begin{longtable}{@{}p{43mm}p{108mm}@{}}
\caption{Glossary.}\label{tab:glossary}\\
\toprule
Term & Plain meaning\\
\midrule\endfirsthead
\toprule
Term & Plain meaning\\
\midrule\endhead
\bottomrule\endfoot
Scalar It\^o SDE & A one-dimensional equation $dX=f(X)\,dt+\sigma(X)\,dW$: a drift $f$ plus noise of size $\sigma$.\\
Locally Lipschitz & On every bounded interval, $f$ and $\sigma$ change at most proportionally to the change in $x$. This makes solutions unique.\\
Maximal solution, lifetime $\zeta$ & The solution followed for as long as it exists; $\zeta$ is the time at which it stops existing.\\
Explosion & $\zeta<\infty$: the path reaches $\pm\infty$ in finite time.\\
Feller's test & An integral test on $f$ and $\sigma$ that decides whether explosion happens.\\
Generator $\cL$ & $\cL V=fV'+\tfrac12\sigma^2V''$. By It\^o's formula, $\cL V$ is the average rate of change of $V(X_t)$.\\
Diffusion zero & A point where $\sigma=0$: the noise switches off there.\\
Degenerate noise & $\sigma$ has zeros; the opposite of nondegenerate noise, $\sigma^2\ge c>0$.\\
Equilibrium & A point with $f=\sigma=0$. A path started there stays there; any other path never hits it.\\
Gate (shunt point) & A point with $\sigma=0$ and $f\neq0$. Paths cross it at most once, in the direction of $f$.\\
Accumulation point, block & A point where infinitely many zeros of $\sigma$ pile up; a block is an interval containing such a pile.\\
Component & A maximal open interval on which $\sigma\neq0$.\\
Scale function $S$ & A change of variable that removes the drift. Hitting probabilities are ratios of its values, as in~\eqref{eq:two-barrier}.\\
Attracting / repelling end & The scale function is finite / infinite at that end. Paths can converge to an attracting end; a repelling end pushes them back.\\
Speed density, speed mass & How long the path lingers near each level; the speed mass is its total integral over a component.\\
Stable in probability & Starting close enough to the equilibrium, the path stays close forever with probability near one.\\
Asymptotically stable & Stable, and the probability of converging to the equilibrium tends to one as the start approaches it.\\
Recurrent & The path returns to every level of its interval infinitely often.\\
Positive / null recurrent & Recurrent with finite / infinite speed mass. Only the positive case has a stationary probability.\\
Stationary (invariant) probability & A distribution of the initial state that the dynamics leaves unchanged.\\
Terminal component & A component that the path eventually never leaves.\\
Selection time & For finitely many gates, the time of the last gate crossing, or zero if no gate is crossed.\\
Weak convergence & Expectations of every bounded continuous function converge.\\
Total-variation convergence & Probabilities of all events converge uniformly; stronger than weak convergence.\\
Empirical measure & The fraction of time up to $t$ that the path spends in each set.\\
Two-point compactification & $[l,r]$ with $\pm\infty$ added as points, so that mass escaping to infinity can be recorded as a point mass there.\\
Green function & The kernel that gives expected time spent near each level before exit; it yields mean exit times.\\
\end{longtable}
}

\section{Local building blocks}\label{sec:local}

Theorem~\ref{thm:A} is assembled from local facts: what the path does near an
equilibrium, near an infinite endpoint and near a gate. This section proves
them. It then describes the path inside a single component, which is the step
that Sections~\ref{sec:finite} and~\ref{sec:accumulation} repeat across many
components. Together the results prove Proposition~\ref{prop:local}.

\subsection{Equilibria}\label{sec:equilibria}

Let $e$ be an equilibrium that is an end of a component $J$, and let
$s\in\{+1,-1\}$ be the side of $e$ on which $J$ lies. Choose $h>0$ with
$e+sy\in J$ for $0<y\le h$. In the \emph{inward coordinate} $Y=s(X-e)$ the
equation becomes
\begin{equation}\label{eq:inward}
 dY_t=\mu(Y_t)\,dt+\nu(Y_t)\,dB_t,\qquad
 \mu(y)=s\,f(e+sy),\quad \nu(y)=\sigma(e+sy),\quad B=sW,
\end{equation}
so that $\mu>0$ means the drift points away from $e$, into $J$. Put
\begin{equation}\label{eq:local-scale}
 \rho_e(y)=\exp\Big(\int_y^h\frac{2\mu(z)}{\nu^2(z)}\,dz\Big),\qquad
 m_e(y)=\frac{2}{\nu^2(y)\rho_e(y)},\qquad 0<y\le h .
\end{equation}
Up to positive constants these are the scale and speed densities of $J$ written
in the coordinate $y$. Thus $e$ is an attracting end of $J$ if and only if
$\int_0^h\rho_e<\infty$, and the speed mass near $e$ is finite if and only if
$\int_0^hm_e<\infty$. We write $H_b$ for the first time $Y$ reaches $b$.

\begin{lemma}[Scale criterion]\label{lem:scale}
The following are equivalent: \textup{(i)} $e$ is an attracting end of $J$;
\textup{(ii)} $e$ is stable in probability from side $s$;
\textup{(iii)} $e$ is asymptotically stable in probability from side $s$.
More precisely, for $0<y<b<h$,
\begin{equation}\label{eq:escape}
 \PP_{e+sy}(H_b<\infty)=
 \begin{cases}
 \displaystyle\int_0^y\rho_e\Big/\int_0^b\rho_e, & \text{if } e \text{ is attracting},\\[8pt]
 1, & \text{if } e \text{ is repelling}.
 \end{cases}
\end{equation}
\end{lemma}

\begin{proof}
For $0<a<y<b<h$ the coefficients are bounded on $[a,b]$ and $\nu$ is bounded
away from zero there. The two-barrier formula~\eqref{eq:two-barrier} in the
coordinate $y$ gives $\PP_{e+sy}(H_b<H_a)=\int_a^y\rho_e/\int_a^b\rho_e$. By
Proposition~\ref{prop:inaccessible}, $Y$ stays positive, so
$\{H_b<\infty\}=\bigcup_n\{H_b<H_{a_n}\}$ for any $a_n\downarrow0$. Letting
$a\downarrow0$ proves~\eqref{eq:escape}.

If $e$ is attracting, the right side of~\eqref{eq:escape} tends to zero as
$y\downarrow0$. A path that moves a distance $\varepsilon>b$ away from $e$ must
first hit $b$, so $e$ is stable. To prove attraction, consider the bounded
martingale $M_t=S(Y_{t\wedge H_b})$, where $S(y)=\int_0^y\rho_e$. It converges
almost surely, with finite quadratic variation
$\langle M\rangle_\infty=\int_0^{H_b}\rho_e^2\nu^2(Y_t)\,dt$. On
$\{H_b=\infty\}$, the limit of $Y_t$ exists in $[0,b]$. A positive limit is
impossible, since $\rho_e^2\nu^2$ is bounded away from zero near it and would
force $\langle M\rangle_\infty=\infty$. Hence $Y_t\to0$ on $\{H_b=\infty\}$, and
$\PP_{e+sy}(X_t\to e)\ge1-\PP_{e+sy}(H_b<\infty)\to1$.

If $e$ is repelling, every nearby initial state reaches the fixed level $b$
almost surely, so $e$ is not stable.
\end{proof}

The threshold for a nondegenerate linearization is the familiar exponent
$\lambda=f'(e)-\tfrac12\sigma'(e)^2$.

\begin{corollary}[First-order test]\label{cor:first-order}
Let $f,\sigma$ be $C^1$ near $e$ with $\sigma'(e)\ne0$. If $\lambda<0$, then $e$
is attracting from both sides and its speed mass is infinite on both sides. If
$\lambda>0$, then $e$ is repelling from both sides with finite speed mass. If
$f,\sigma$ are $C^2$ and $\lambda=0$, then $e$ is repelling with infinite speed
mass.
\end{corollary}

\begin{proof}
In~\eqref{eq:inward}, $\mu(y)=f'(e)y+o(y)$ and
$\nu^2(y)=\sigma'(e)^2y^2(1+o(1))$ on both sides, so
$2\mu/\nu^2=k/y+o(1/y)$ with $k=2f'(e)/\sigma'(e)^2$, and $\lambda<0$ if and only
if $k<1$. For each $\varepsilon>0$ this gives
$c\,y^{-(k-\varepsilon)}\le\rho_e(y)\le C\,y^{-(k+\varepsilon)}$ near zero,
and hence $m_e(y)$ lies between multiples of $y^{k\mp\varepsilon-2}$. For
$k\ne1$ both integrability questions are therefore decided by the sign of
$k-1$. Under $C^2$ regularity the error is $O(1)$, so $\rho_e(y)\sim Cy^{-k}$
and $m_e(y)\sim Dy^{k-2}$; at $k=1$ both integrals diverge.
\end{proof}

\begin{example}[Noise-induced stability]\label{ex:noise-stability}
For $dX_t=\tfrac14X_t\,dt+X_t\,dW_t$ the drift points away from zero, yet
$\lambda=-\tfrac14$, so zero is attracting from both sides. A quadratic
Lyapunov function cannot show this, since $V(x)=x^2$ gives $\cL V=\tfrac32x^2>0$.
The function $V(x)=|x|^p$ with $0<p<\tfrac12$ gives
$\cL V=p(\tfrac14+\tfrac{p-1}2)|x|^p<0$ and does (Appendix~\ref{app:lyapunov}).
\end{example}

\begin{example}[$C^1$ does not decide $\lambda=0$]\label{ex:critical}
For $\sigma(x)=x$ and $f(x)=x/2$ we have $\lambda=0$ and $X_t=xe^{W_t}$, so every
initial state reaches any outer level: zero is repelling. Now keep
$\sigma(x)=x$ but take $f(x)=x/2-x/\log(1/|x|)$ for $0<|x|<e^{-4}$, extended to
globally Lipschitz coefficients. The derivatives at zero are unchanged, so still
$\lambda=0$, but
$2\mu/\nu^2=1/y-2/(y\log(1/y))$ and
\[
 \rho_e(y)=C\,\frac{1}{y\,\log^2(1/y)},
\]
which is integrable at zero. Zero is now attracting from both sides.
\end{example}

When the linearization degenerates, the orders of vanishing decide. Write, in
the inward coordinate,
\begin{equation}\label{eq:orders}
 \mu(y)=Ay^p+O(y^{p+1}),\qquad \nu(y)=By^n+O(y^{n+1}),\qquad AB\ne0,
\end{equation}
with integers $p\ge0$, $n\ge1$, and put $k=2A/B^2$.

\begin{theorem}[Finite-order tests]\label{thm:finite-order}
Under~\eqref{eq:orders}, near $y=0$:
\begin{enumerate}[label=\textup{(\roman*)}]
\item if $p<2n-1$, then $\log\rho_e(y)=\frac{k}{q}\,y^{-q}(1+o(1))$ with
$q=2n-1-p>0$. For $A>0$ the end is repelling with finite speed mass; for $A<0$
it is attracting with infinite speed mass;
\item if $p=2n-1$, then $\rho_e(y)\sim Cy^{-k}$ and $m_e(y)\sim Dy^{k-2n}$. The
end is repelling if and only if $k\ge1$, and its speed mass is finite if and
only if $k>2n-1$;
\item if $p>2n-1$, then $\rho_e$ has a positive finite limit and
$m_e(y)\sim Dy^{-2n}$. The end is attracting with infinite speed mass.
\end{enumerate}
Part~\textup{(iii)} also holds under the diffusion expansion in~\eqref{eq:orders}
and the drift bound $\mu(y)=O(y^{2n})$, including $\mu\equiv0$.
\end{theorem}

\begin{proof}
By~\eqref{eq:orders}, $2\mu/\nu^2=ky^{p-2n}(1+O(y))$. If $p<2n-1$, integrating
gives the stated asymptotics of $\log\rho_e$, with an error of smaller order
than $y^{-q}$. For $A>0$, $\rho_e\ge e^{cy^{-q}}$ is not integrable, while
$m_e\le Cy^{-2n}e^{-cy^{-q}}$ is. For $A<0$ the inequalities reverse. If
$p=2n-1$, the ratio is $k/y+O(1)$, and integrating gives $\rho_e\sim Cy^{-k}$;
the power criterion decides both integrals. If $p>2n-1$, the ratio is bounded,
so $\rho_e$ has a positive limit. The same boundedness follows directly from
$\mu(y)=O(y^{2n})$ and $\nu(y)\sim By^n$.
\end{proof}

At an equilibrium with $f(e+z)=a_mz^m+O(z^{m+1})$ and
$\sigma(e+z)=b_nz^n+O(z^{n+1})$, the inward coefficients are $A=s^{m+1}a_m$ and
$B^2=b_n^2$, with $p=m$. In particular, in the critical case $m=2n-1$ the
exponent $m$ is odd and both sides behave alike. The two thresholds in~(ii)
answer different questions: for $n\ge2$ and $1\le k\le2n-1$ the end is
repelling, yet its speed mass is infinite. This is the source of null recurrence
in Theorem~\ref{thm:B}. No finite order decides a flat zero: with
$\sigma(x)=e^{-1/x^2}$ and $f=\pm x\sigma$, all derivatives agree at zero, yet
the minus sign gives an attracting and the plus sign a repelling equilibrium.

\begin{corollary}[Boundary profile of a stationary density]\label{cor:profile}
Let $J$ be positive recurrent and $e$ an end of $J$. In case~\textup{(i)} of
Theorem~\ref{thm:finite-order} with $A>0$, the stationary density satisfies
$\log\pi_J(e+sy)=-\frac kq\,y^{-q}(1+o(1))$, so it vanishes faster than any
power. In case~\textup{(ii)} with $k>2n-1$, $\pi_J(e+sy)\sim C_ey^{k-2n}$: it is
unbounded if $k<2n$, has a positive limit if $k=2n$, and vanishes if $k>2n$.
\end{corollary}

\begin{proof}
The stationary density is proportional to $m_e=2/(\nu^2\rho_e)$, and
$\log\nu^2=O(\log y)$ is negligible against $y^{-q}$.
\end{proof}

\begin{example}[Both thresholds in one family]\label{ex:beta}
For an integer $n\ge1$ and $\kappa>0$, consider on $(0,1)$
\[
 dX_t=\tfrac{\kappa}{2}X_t^{2n-1}(1-X_t)^{2n-1}(1-2X_t)\,dt
 +X_t^n(1-X_t)^n\,dW_t .
\]
Here $2f/\sigma^2=\kappa\big(\frac1x-\frac1{1-x}\big)$, so
$\rho\propto[x(1-x)]^{-\kappa}$ and $m\propto[x(1-x)]^{\kappa-2n}$. The process
is recurrent if and only if $\kappa\ge1$. It is null recurrent for
$1\le\kappa\le2n-1$ and positive recurrent for $\kappa>2n-1$, with stationary
density the $\mathrm{Beta}(\alpha,\alpha)$ density, $\alpha=\kappa-2n+1$. For
$n=2$ the density is unbounded at both ends for $3<\kappa<4$, uniform for
$\kappa=4$, and vanishes at both ends for $\kappa>4$. For $0<\kappa<1$ both ends
attract, and Theorem~\ref{thm:component} gives convergence to a random endpoint.
\end{example}

\subsection{Infinite endpoints}\label{sec:infinity}

At $+\infty$ the same two questions, attraction and speed mass, are decided by
growth rates. Suppose that, as $x\to\infty$,
\begin{equation}\label{eq:growth}
 f(x)=Ax^{p}\big(1+O(x^{-\delta})\big),\qquad
 \sigma^2(x)=B^2x^{2n}\big(1+O(x^{-\delta})\big),
\end{equation}
with real exponents $p$ and $n\ge0$, constants $A\ne0$, $B\ne0$ and some
$\delta>0$. Put $k=2A/B^2$ and $q=p-2n+1$.

\begin{proposition}[Tests at infinity]\label{prop:infinity}
Let $+\infty$ be the right end of a component $J$ and assume~\eqref{eq:growth}.
\begin{enumerate}[label=\textup{(\roman*)}]
\item If $q>0$, then $+\infty$ is attracting with infinite speed mass when
$A>0$, and repelling with finite speed mass when $A<0$.
\item If $q=0$, then $\rho_J(x)\sim Cx^{-k}$ and $m_J(x)\sim Dx^{k-2n}$:
$+\infty$ is attracting if and only if $k>1$, and its speed mass is finite if and
only if $k<2n-1$.
\item If $q<0$, then $\rho_J$ has a positive finite limit, so $+\infty$ is
repelling, and $m_J(x)\sim Dx^{-2n}$ has finite mass near $+\infty$ if and only
if $n>\tfrac12$.
\end{enumerate}
The end $-\infty$ is treated by the reflection $x\mapsto-x$.
\end{proposition}

\begin{proof}
By~\eqref{eq:growth}, $2f/\sigma^2=kx^{p-2n}(1+O(x^{-\delta}))$. If $q>0$,
$\int_c^x2f/\sigma^2=\frac kqx^{q}(1+o(1))$; for $A>0$ the scale density decays
like $e^{-cx^q}$ and the speed density grows like $x^{-2n}e^{cx^q}$, and for
$A<0$ the roles reverse. If $q=0$, the ratio is $k/x+O(x^{-1-\delta})$, whose
error is integrable, so $\rho_J\sim Cx^{-k}$. If $q<0$, the ratio itself is
integrable at infinity.
\end{proof}

For instance, the stochastic logistic equation
$dX=X(\alpha-\beta X)\,dt+\tau X\,dW$ has $p=2$, $n=1$ and $A=-\beta$, so
$+\infty$ is repelling with finite speed mass. When $+\infty$ is attracting,
the path escapes to infinity with positive probability, but it may or may not
do so in finite time. By Feller's test~\cite[Ch.~15]{Karlin},
\cite[\S5.5]{KaratzasShreve}, explosion toward $+\infty$ has positive
probability if and only if
\begin{equation}\label{eq:feller}
 \int_c^\infty\big(S_J(\infty)-S_J(y)\big)\,m_J(y)\,dy<\infty .
\end{equation}
In particular, coefficients of at most linear growth never explode.

\subsection{Gates}\label{sec:gates}

A gate behaves very differently from an equilibrium. It is visited, and once
crossed it is never crossed back.

\begin{theorem}[One-way passage]\label{thm:one-way}
Let $f,\sigma$ be locally Lipschitz and let $h$ be a gate, $\sigma(h)=0$,
$f(h)\neq0$. For the maximal solution started at $x$:
\begin{enumerate}[label=\textup{(\roman*)}]
\item if $f(h)>0$ and $x\ge h$, then $X_t>h$ for all $0<t<\zeta$, almost surely;
if $f(h)<0$ and $x\le h$, then $X_t<h$ for all $0<t<\zeta$;
\item on $\{T_h<\zeta\}$, $(X_{T_h+t}-h)/t\to f(h)$ as $t\downarrow0$, and
$X_{T_h+t}-h$ has the sign of $f(h)$ for all $0<t<\zeta-T_h$; in particular $h$
is visited at most once;
\item if $f\ge\beta>0$ on $[h-\delta,h]$, then for $h-\delta<x<h$,
\[
 \EE_x\big(T_h\wedge T_{h-\delta}\big)\le\frac{h-x}{\beta},\qquad
 \PP_x\big(T_h<T_{h-\delta}\big)\ge\frac{x-h+\delta}{\delta}.
\]
The case $f(h)<0$ is symmetric.
\end{enumerate}
\end{theorem}

\begin{proof}
Put $b=f(h)$, $Y=X-h$, and
$a(y)=(f(h+y)-b)/y$, $c(y)=\sigma(h+y)/y$ for $y\ne0$, both locally bounded.
Then $dY=(b+a(Y)Y)\,dt+c(Y)Y\,dW$. With the positive process
$E_t=\exp\{\int_0^t(a-\tfrac12c^2)(Y_s)\,ds+\int_0^tc(Y_s)\,dW_s\}$, It\^o's
product rule gives $d(E_t^{-1}Y_t)=bE_t^{-1}\,dt$, hence
\begin{equation}\label{eq:variation}
 X_t-h=E_t\Big(x-h+b\int_0^tE_s^{-1}\,ds\Big),\qquad 0\le t<\zeta,
\end{equation}
simultaneously for all $t$, almost surely. The bracket is strictly monotone with
the sign of $b$, which proves~(i) and shows it vanishes at most once. At
$T_h$ it vanishes, and~\eqref{eq:variation} restarted there gives~(ii). For~(iii),
stop at $\tau=T_h\wedge T_{h-\delta}$; the stopped stochastic integral is a
martingale, so $h\ge\EE_xX_{t\wedge\tau}\ge x+\beta\,\EE_x(t\wedge\tau)$.
Letting $t\to\infty$ gives the mean bound, and $\EE_xX_\tau\ge x$ with
$X_\tau\in\{h-\delta,h\}$ gives the probability bound.
\end{proof}

\begin{example}[An explicit gate]\label{ex:explicit-gate}
For $dX_t=dt+X_t\,dW_t$,
\[
 X_t=e^{W_t-t/2}\Big(x+\int_0^te^{-W_s+s/2}\,ds\Big).
\]
For $x\ge0$ the path is positive at all positive times, and from $x=0$,
$X_t/t\to1$. For $x<0$ the integral grows without bound, since $W_s/s\to0$, so the
path hits zero at a unique finite time and stays positive afterwards.
\end{example}

The next lemma gives the gate rows of Table~\ref{tab:local}. It uses only the
sign of $f$ at the gate and the Lipschitz bound $|\sigma(z)|\le L|z-h|$.

\begin{lemma}[Gate ends]\label{lem:gate-ends}
Let $v$ be a gate that is the right end of a component $J$, and let $S=S_J$,
$m=m_J$.
\begin{enumerate}[label=\textup{(\roman*)}]
\item If $f(v)<0$ (drift into $J$), then $S(v-)=+\infty$ and
$\int^vm<\infty$: the end is repelling with finite speed mass.
\item If $f(v)>0$ (drift out of $J$), then $S(v-)<\infty$ and
$\int^v\big(S(v-)-S(y)\big)\,m(y)\,dy<\infty$.
\end{enumerate}
The same holds at a left end with the signs of $f$ reversed.
\end{lemma}

\begin{proof}
Choose $\delta,\beta,L>0$ with $|f|\ge\beta$ and $|\sigma(z)|\le L(v-z)$ on
$[v-\delta,v)$. In case~(i), for $v-\delta\le y<v$,
\[
 \log\frac{\rho(y)}{\rho(v-\delta)}=-\int_{v-\delta}^y\frac{2f}{\sigma^2}
 \ge\frac{2\beta}{L^2}\Big(\frac1{v-y}-\frac1\delta\Big),
\]
so $\rho\to\infty$ and $S(v-)=\infty$. Let $g=1/\rho$; then $g'=fm$, so
$\beta m\le-g'$ and $\int_{v-\delta}^vm\le g(v-\delta)/\beta<\infty$. In
case~(ii) the same estimate with reversed sign gives
$\rho(y)\le Ce^{-c/(v-y)}$, so $S(v-)<\infty$. With $R(y)=S(v-)-S(y)$ and
$g=1/\rho$ we have $R\,g\le v-y\to0$ because $\rho$ decreases near $v$, and
$g'=fm\ge\beta m$. Integrating $R\,g'$ by parts gives
$\beta\int_{v-\delta}^vRm\le\delta$.
\end{proof}

An attracting gate end is reached in finite time, not merely approached.

\begin{lemma}[Attracting gates are hit]\label{lem:gate-hit}
Let $v$ be a gate end of a component $J$ with drift out of $J$, and let $\tau_J$
be the exit time from $J$. Then, almost surely,
$\{\lim_{t\uparrow\tau_J}X_t=v\}\subseteq\{\tau_J<\infty,\ X_{\tau_J}=v\}$.
\end{lemma}

\begin{proof}
Consider a right end $v$ and choose the drift strip $(v-\delta,v)$ as in
Theorem~\ref{thm:one-way}(iii). Starting anywhere in this strip, its exit time
has finite mean. On $\{\tau_J=\infty,\ X_t\to v\}$ the path stays in the strip
after some integer time. For each fixed integer the Markov property gives
probability zero to this event of permanent residence. Taking their countable
union excludes infinite-time approach. A finite approach to the interior
point $v$ gives $X_{\tau_J}=v$ by continuity and local continuation. The
left-end proof follows by reflection.
\end{proof}

\subsection{One component}\label{sec:component}

We now combine the three kinds of ends. Let $J=(u,v)$ be a component, so that
each end is an equilibrium, an infinite endpoint or a gate. Let
$\tau_J=\inf\{t<\zeta:X_t\notin J\}$, with $\tau_J=\zeta$ if the path never
leaves $J$, and let $S=S_J$.

\begin{theorem}[Behavior in one component]\label{thm:component}
Let $x\in J$.
\begin{enumerate}[label=\textup{(\roman*)}]
\item If both ends are repelling, then the path never leaves $J$ and
$\zeta=\infty$. It is recurrent, with $\liminf X_t=u$ and $\limsup X_t=v$, and it
is positive recurrent if and only if $M_J<\infty$
\textup{(}Lemma~\ref{lem:recurrence}\textup{)}.
\item If $u$ is repelling and $v$ attracting, then $X_t\to v$ as
$t\uparrow\tau_J$ almost surely, and for $u<a<x$,
\begin{equation}\label{eq:mixed}
 \PP_x(T_a<\tau_J)=\frac{S(v-)-S(x)}{S(v-)-S(a)} .
\end{equation}
\item If both ends are attracting, then $X_t$ converges to $u$ or $v$ as
$t\uparrow\tau_J$, and
\begin{equation}\label{eq:split}
 \PP_x\Big(\lim_{t\uparrow\tau_J}X_t=v\Big)=\frac{S(x)-S(u+)}{S(v-)-S(u+)}\in(0,1).
\end{equation}
\end{enumerate}
In~\textup{(ii)} and~\textup{(iii)}, on the event that a given end is the
limit, a gate is hit in finite time and crossed for good, a finite
equilibrium is approached at infinite time, and an infinite endpoint is
approached as $t\uparrow\zeta$. Explosion toward an infinite endpoint has
positive probability exactly when its Feller integral~\eqref{eq:feller}
(or its reflection) is finite. The case of $u$ attracting and $v$ repelling
is symmetric.
\end{theorem}

\begin{proof}
By Itô's formula, $S(X_t)$ is a continuous local martingale on $[0,\tau_J)$.

(i) A repelling equilibrium end is never reached
(Proposition~\ref{prop:inaccessible}). A repelling gate end has drift into $J$
(Lemma~\ref{lem:gate-ends}) and is never crossed (Theorem~\ref{thm:one-way}).
A repelling infinite end is not reached in finite time, since~\eqref{eq:feller}
needs $S(\infty)<\infty$. Hence the process stays in $J$ and is conservative, and
Lemma~\ref{lem:recurrence} applies. Letting one barrier tend to an end
in~\eqref{eq:two-barrier} shows that every level is hit, and the alternating
visits to levels approaching $u$ and $v$ give the $\liminf$ and $\limsup$.

(iii) Now $S(X_t)$ is bounded on $[0,\tau_J)$, so it converges almost surely and
its quadratic variation $\int_0^{\tau_J}(\rho\sigma)^2(X_t)\,dt$ is finite. The
inverse of $S$ gives a limit in the extended interval $[u,v]$. An interior
limit at finite $\tau_J$ would allow local continuation within $J$. At
infinite $\tau_J$, an interior limit would force infinite quadratic variation,
since $(\rho\sigma)^2$ is bounded below near that limit. Thus only the two ends
are possible. Bounded convergence in scale gives~\eqref{eq:split}.

(ii) Now $Z_t=S(v-)-S(X_t)$ is a nonnegative local martingale and converges
almost surely as $t\uparrow\tau_J$ by nonnegative supermartingale convergence
after localization. For each integer $n>Z_0$, stop upon reaching $n$.
The resulting bounded scale martingale has finite terminal quadratic
variation. The union of the events $\{\sup_{t<\tau_J} Z_t<n\}$ has full
probability, so the original quadratic variation is finite as well.
As in~(iii), the limit is an end, and $u$ is excluded since $Z$ diverges
there. Stopping also at $T_a$ and taking the terminal expectation gives
\eqref{eq:mixed}.

Finally, a gate end that is the limit is hit in finite time by
Lemma~\ref{lem:gate-hit}.
\end{proof}

\begin{proposition}[Mean exit through gates]\label{prop:mean-exit}
Let $J=(u,v)$ be bounded, with both ends gates and at least one of them
attracting. Then $\sup_{x\in J}\EE_x\tau_J<\infty$, and
$\EE_x\tau_J=\int_JG_J(x,y)\,m_J(y)\,dy$, where
\[
 G_J(x,y)=
 \begin{cases}
 \dfrac{\big(S(x\wedge y)-S(u+)\big)\big(S(v-)-S(x\vee y)\big)}{S(v-)-S(u+)},
 & \text{both ends attracting},\\[10pt]
 S(v-)-S(x\vee y), & \text{only }v\text{ attracting},\\[3pt]
 S(x\wedge y)-S(u+), & \text{only }u\text{ attracting}.
 \end{cases}
\]
If $u$ is an inward gate and $v$ an outward gate, then
$\EE_uT_v\le\sup_{x\in J}\EE_x\tau_J$. The reflected statement holds for a
start at an inward right gate.
\end{proposition}

\begin{proof}
The Green identity follows by compact exhaustion in Lemma~\ref{lem:green}.
When both gates are outward, the kernel is bounded above by $S(y)-S(u+)$
near $u$ and by $S(v-)-S(y)$ near $v$, uniformly in $x$.
Lemma~\ref{lem:gate-ends}(ii) makes both bounds integrable against $m_J$.

Suppose $u$ points inward and $v$ outward. Put $B=S(v-)$ and $g=1/\rho_J$.
Choose $c\in J$ and $\beta>0$ with $f\ge\beta$ on $[u,c]$.
By the inward-gate estimate, $g(u+)=0$. Since $g'=fm_J$,
$\int_u^z m_J(y)\,dy\le g(z)/\beta$ for $u<z\le c$. Tonelli's theorem gives
\begin{align*}
 \int_u^c (B-S(y))m_J(y)\,dy
 &=(B-S(c))\int_u^c m_J(y)\,dy+\int_u^c \rho_J(z)\left(\int_u^z m_J(y)\,dy\right)dz .
\end{align*}
The second term is at most $(c-u)/\beta$, and the first is finite.
Lemma~\ref{lem:gate-ends}(ii) gives integrability near $v$. The mixed-case
supremum in Lemma~\ref{lem:green} is therefore finite.

From $u$, Theorem~\ref{thm:one-way} gives $X_\varepsilon\in J$ on
$\{T_v>\varepsilon\}$. Conditioning at $\varepsilon$ yields
$\EE_uT_v\le\varepsilon+\sup_{x\in J}\EE_x\tau_J$. Let $\varepsilon\downarrow0$.
Reflection treats an inward right gate.
\end{proof}

With a repelling equilibrium at one end and an outward gate at the other,
Theorem~\ref{thm:component}(ii) still gives almost-sure passage to the gate.
Section~\ref{sec:times} determines when its mean is finite.

\begin{proof}[Proof of Proposition~\ref{prop:local}]
The first row of Table~\ref{tab:local} is Corollary~\ref{cor:first-order}. The
other equilibrium rows are Theorem~\ref{thm:finite-order}, the gate rows are
Lemma~\ref{lem:gate-ends} with Lemma~\ref{lem:gate-hit}, and the rows at
infinity are Proposition~\ref{prop:infinity}.
\end{proof}

The following examples illustrate the configurations of
Theorem~\ref{thm:component}.

\begin{example}[Two attracting equilibria, no drift]\label{ex:zero-drift}
For $dX_t=X_t(1-X_t)\,dW_t$ on $(0,1)$, the scale is $S(x)=x$. Both ends
attract, and $X_t\to X_\infty\in\{0,1\}$ with $\PP_x(X_\infty=1)=x$. The path
stays strictly inside $(0,1)$ at every finite time.
\end{example}

\begin{example}[One attracting equilibrium]\label{ex:one-attracting}
For $dX_t=\kappa X_t(1-X_t)\,dt+X_t(1-X_t)\,dW_t$ on $(0,1)$ with
$\kappa\ge\tfrac12$, $\rho(x)\propto x^{-2\kappa}(1-x)^{2\kappa}$. It is not
integrable at zero and is integrable at one, so $X_t\to1$ almost surely, including
at the critical value $\kappa=\tfrac12$.
\end{example}

\begin{example}[Two inward gates]\label{ex:inward}
For $dX_t=(1-2X_t)\,dt+X_t(1-X_t)\,dW_t$, the points $0$ and $1$ are gates with
drift into $(0,1)$, and $\rho(x)\propto\exp\{2/(x(1-x))\}$. The process is
positive recurrent on $(0,1)$ with stationary density
\[
 \pi(x)=\frac{C}{x^2(1-x)^2}\exp\Big(-\frac2{x(1-x)}\Big),
\]
and $\liminf X_t=0$, $\limsup X_t=1$ almost surely.
\end{example}

\begin{example}[Two outward gates]\label{ex:outward}
For $dX_t=(2X_t-1)\,dt+X_t(1-X_t)\,dW_t$, the gates $0$ and $1$ have drift out of
$(0,1)$, and $\rho(x)\propto\exp\{-2/(x(1-x))\}$. The path leaves $(0,1)$ in a
time of finite mean, through $1$ with probability
$\int_0^x\rho\big/\int_0^1\rho$, which equals $\tfrac12$ at $x=\tfrac12$.
It then stays on the far side of the gate it crossed.
\end{example}

\begin{example}[An explicit mean exit time]\label{ex:mean-exit}
For $dX_t=dt+X_t(1-X_t)\,dW_t$ on $(0,1)$, zero is an inward gate and one is an
outward gate. Since the stopped drift identity gives
$\EE_xX_{t\wedge T_1}=x+\EE_x(t\wedge T_1)$, we get $\EE_xT_1=1-x$.
\end{example}

\section{Finitely many gates: directed routes}\label{sec:finite}

When $Z$ is finite, Theorem~\ref{thm:component} applies component by component,
and the path follows a route that can be read off a directed graph. This section
describes that route. It is the finite case of Theorem~\ref{thm:A}, and it shows
what the accumulation argument of Section~\ref{sec:accumulation} must replace.

Assume~\eqref{eq:standing} and $Z=\{h_1<\dots<h_N\}$. Put $h_0=l$, $h_{N+1}=r$,
$I_i=(h_i,h_{i+1})$ for $0\le i\le N$, and $\varepsilon_j=\operatorname{sgn}f(h_j)$.

\begin{definition}[Route graph]\label{def:route}
The \emph{route graph} has the components $I_0,\dots,I_N$ as vertices, and one
edge at each gate: $I_{j-1}\to I_j$ if $\varepsilon_j=+1$ and $I_j\to I_{j-1}$ if
$\varepsilon_j=-1$. An \emph{exit} of $I_i$ is an end through which the path can
leave or converge out of it: a gate end with drift out of $I_i$, or an outer end
($l$ for $I_0$, $r$ for $I_N$) that is attracting. A component is
\emph{terminal} if it has no exit.
\end{definition}

\begin{figure}[htbp]
\centering
\resizebox{0.98\linewidth}{!}{%
\begin{tikzpicture}[>=Stealth,font=\small,
 cell/.style={draw=black!35,fill=black!3,rounded corners=2mm,
   minimum height=1.1cm,align=center},
 terminal/.style={cell,draw=teal!70!black,fill=teal!9,line width=0.9pt},
 gate/.style={->,line width=1pt,draw=black!75}]
\node[cell,text width=1.75cm] at (0,0) {$I_0$\\[-1pt]{\scriptsize transient}};
\node[terminal,text width=1.75cm] at (2.75,0) {$I_1$\\[-1pt]{\scriptsize positive recurrent}};
\node[cell,text width=1.75cm] at (5.5,0) {$I_2$\\[-1pt]{\scriptsize two exits}};
\node[cell,text width=1.75cm] at (8.25,0) {$I_3$\\[-1pt]{\scriptsize transient}};
\node[terminal,text width=1.75cm] at (11,0) {$I_4$\\[-1pt]{\scriptsize positive recurrent}};
\node[cell,text width=1.75cm] at (13.75,0) {$I_5$\\[-1pt]{\scriptsize transient}};
\draw[gate] (1.02,0) -- (1.73,0);
\node[font=\scriptsize] at (1.375,0.45) {$+$};
\node[font=\scriptsize] at (1.375,-0.45) {$h_1$};
\draw[gate] (4.48,0) -- (3.77,0);
\node[font=\scriptsize] at (4.125,0.45) {$-$};
\node[font=\scriptsize] at (4.125,-0.45) {$h_2$};
\draw[gate] (6.52,0) -- (7.23,0);
\node[font=\scriptsize] at (6.875,0.45) {$+$};
\node[font=\scriptsize] at (6.875,-0.45) {$h_3$};
\draw[gate] (9.27,0) -- (9.98,0);
\node[font=\scriptsize] at (9.625,0.45) {$+$};
\node[font=\scriptsize] at (9.625,-0.45) {$h_4$};
\draw[gate] (12.73,0) -- (12.02,0);
\node[font=\scriptsize] at (12.375,0.45) {$-$};
\node[font=\scriptsize] at (12.375,-0.45) {$h_5$};
\node[font=\scriptsize,align=center] at (6.875,-1.35)
 {$I_i=(h_i,h_{i+1}),\quad h_0=l,\ h_6=r;\qquad +\text{ or }-\ =\ \operatorname{sgn}f(h_j)$};
\node[font=\small,align=center] at (6.875,-2.05)
 {Starting at $x\in I_2$: one random branch, then one-way passage.};
\node[terminal,text width=5.0cm] at (2.3,-3.2)
 {Selected interval $I_1$\\[-1pt]{\scriptsize $w_L(x)=\dfrac{B-S(x)}{B-A}$}};
\node[terminal,text width=5.0cm] at (11.45,-3.2)
 {Selected interval $I_4$\\[-1pt]{\scriptsize $w_R(x)=\dfrac{S(x)-A}{B-A}$}};
\node[font=\scriptsize,align=center] at (6.875,-4.45)
 {$S$ is the scale function of $I_2=(h_2,h_3)$;\quad $A=S(h_2+)$,\ $B=S(h_3-)$.};
\node[font=\scriptsize,align=center] at (6.875,-5.1)
 {Each gate is crossed at most once. Both selected intervals are bounded by two inward gates, hence positive recurrent.};
\end{tikzpicture}}
\caption{A route graph with five gates and signs $(+,-,+,+,-)$, assuming both outer
ends repel. The terminal components are $I_1$ and $I_4$. A path started in $I_2$
makes one random choice, with the scale probabilities shown; every later step is
forced.}
\label{fig:route}
\end{figure}
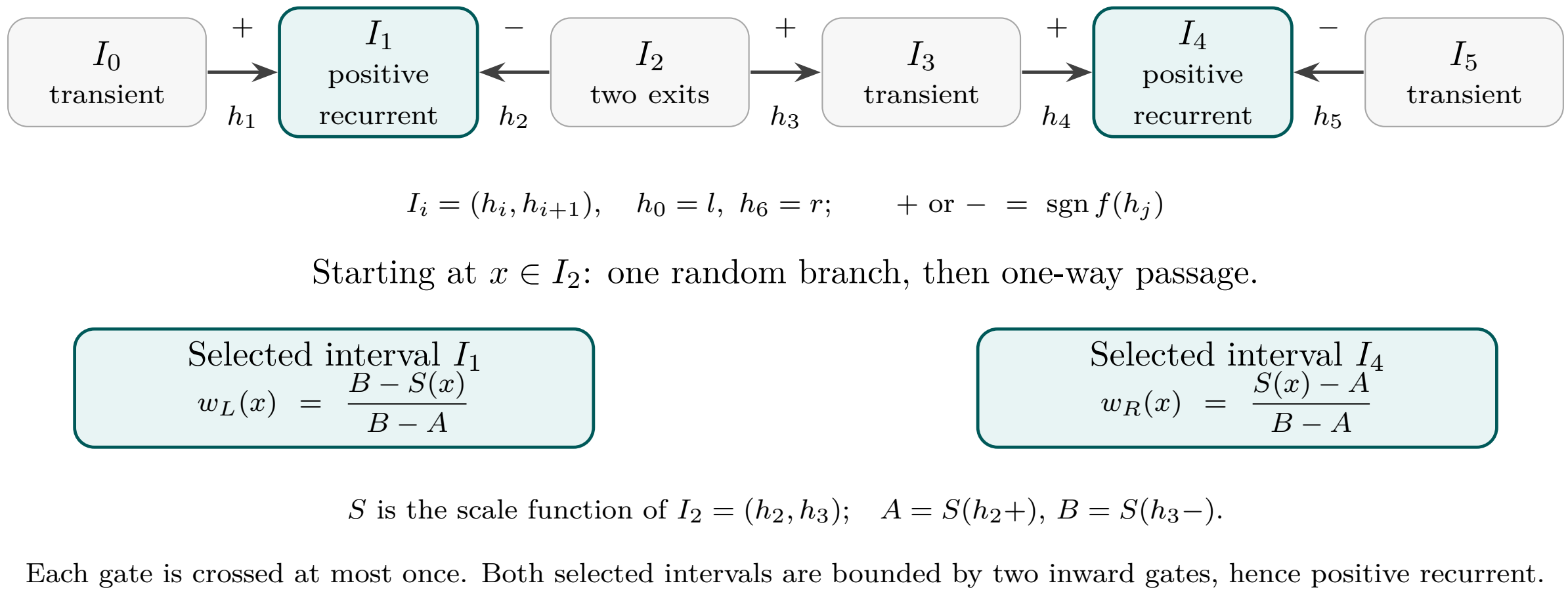

By Lemma~\ref{lem:gate-ends}, the exits of $I_i$ are exactly its attracting ends.
Hence a component is terminal if and only if both its ends repel, that is, if and
only if it belongs to $\Rc$. Figure~\ref{fig:route} shows an example.

\begin{theorem}[Directed routes]\label{thm:routes}
Let $x\in(l,r)$.
\begin{enumerate}[label=\textup{(\roman*)}]
\item The path visits components along a directed path of the route graph. Each
gate is crossed at most once, at a finite time, and at most $N$ gates are
crossed. A path started at a gate $h_j$ enters $I_j$ immediately if
$\varepsilon_j=+1$ and $I_{j-1}$ if $\varepsilon_j=-1$.
\item The route ends in one of two ways: the path stays forever in a terminal
component, where it is recurrent; or it converges to an attracting outer end
as $t\uparrow\zeta$.
\item Only the starting component can have two exits. If its exits are both
ends $u<v$, the path leaves through $v$ with probability
$\big(S(x)-S(u+)\big)/\big(S(v-)-S(u+)\big)$, computed with the scale function
of the starting component; every later step is deterministic. Hence each
initial state has at most two possible final outcomes.
\item A terminal component bounded by two gates is positive recurrent. If
$N\ge1$, a null recurrent terminal component is $I_0$ or $I_N$, and its
speed mass is infinite only at its outer end. If $N=0$, a terminal component
is null recurrent exactly when $M_{(l,r)}=\infty$; either or both ends may
carry infinite speed mass.
\end{enumerate}
\end{theorem}

\begin{proof}
Apply Theorem~\ref{thm:component} in the starting component. If it is terminal,
the path stays there. Otherwise the path either converges to an attracting outer
end or crosses a gate end at a finite time. In the second case, by
Theorem~\ref{thm:one-way}, it enters the adjacent component and never returns.
The gate it crossed is an end of the new component with drift into it, so the new
component has at most one exit, its other end. Theorem~\ref{thm:component}(i)
or~(ii) applies again, and the strong Markov property at the crossing time
continues the argument. The route graph is a path with one orientation per edge,
so no component is entered twice, and at most $N$ crossings occur. This
proves~(i)--(iii); the probability in~(iii) is~\eqref{eq:split}. For a start at a
gate, apply Theorem~\ref{thm:one-way}(ii) and the Markov property at small
positive times.

For~(iv), the speed mass of a component is finite near a gate end with inward
drift (Lemma~\ref{lem:gate-ends}) and on compact subintervals, by continuity of
$m_J$. A component bounded by two inward gates therefore has $M_J<\infty$.
When $N\ge1$, each remaining terminal component has just one outer end.
When $N=0$, Lemma~\ref{lem:recurrence} applies on $(l,r)$.
\end{proof}

Theorem~\ref{thm:routes} contains the cases of the original
classification, in which both outer ends repel and there are at most two gates.
Table~\ref{tab:small} lists them.

\begin{corollary}[Zero, one and two gates]\label{cor:small}
Assume both outer ends are repelling. Then the terminal components are as in
Table~\ref{tab:small}. In the last row, a path started in $(h_1,h_2)$ selects
$(h_2,r)$ with probability $w_R(x)=(S(x)-S(h_1+))/(S(h_2-)-S(h_1+))$, computed
with the scale function of $(h_1,h_2)$, and $(l,h_1)$ with probability
$1-w_R(x)$.
\end{corollary}

\begin{table}[htbp]
\centering\small
\renewcommand{\arraystretch}{1.15}
\begin{tabular}{@{}lll@{}}
\toprule
Gates and signs & Terminal component(s) & Recurrence type\\
\midrule
none & $(l,r)$ & positive iff $M_{(l,r)}<\infty$\\
$h$, $+$ & $(h,r)$ & positive iff the speed mass is finite near $r$\\
$h$, $-$ & $(l,h)$ & positive iff the speed mass is finite near $l$\\
$h_1,h_2$, $(+,-)$ & $(h_1,h_2)$ & positive\\
$h_1,h_2$, $(+,+)$ & $(h_2,r)$ & positive iff the speed mass is finite near $r$\\
$h_1,h_2$, $(-,-)$ & $(l,h_1)$ & positive iff the speed mass is finite near $l$\\
$h_1,h_2$, $(-,+)$ & $(l,h_1)$ and $(h_2,r)$ & each by the speed mass at its outer end\\
\bottomrule
\end{tabular}
\caption{Terminal components with at most two gates and repelling outer ends
(Corollary~\ref{cor:small}). In every row except the last, the terminal
component is selected almost surely from every initial state.}
\label{tab:small}
\end{table}

The last example shows that a point where $\sigma\neq0$ can be locally regular
yet unreachable from much of the state space, and it exercises the infinite
endpoints.

\begin{example}[A one-way barrier on the whole line]\label{ex:barrier}
Let $f\equiv1$ and
\[
 \sigma(u)=\frac{u(1-u^2)}{1+u^2},
\]
which is smooth and globally Lipschitz. There is no equilibrium, so $l=-\infty$
and $r=+\infty$, and the gates are $-1$, $0$ and $1$, all rightward. The point
$c=-\tfrac12$ is regular inside $(-1,0)$, yet a path started at $x=\tfrac12$ stays
positive (Theorem~\ref{thm:one-way}), so $\PP_{1/2}(T_{-1/2}<\infty)=0$. As
$|u|\to\infty$, $\sigma(u)\sim-u$ and $f=1$, so $p=0<1=2n-1$ in
Proposition~\ref{prop:infinity}: both infinite ends repel, and the speed density
behaves like $u^{-2}$. By Theorem~\ref{thm:routes}, every path crosses the
remaining gates to its right and ends in $(1,\infty)$, where it is positive
recurrent.
\end{example}

\section{Accumulating zeros: the global alternative}\label{sec:accumulation}

Theorem~\ref{thm:routes} counts gates: each is crossed at most once, so after
finitely many crossings the path settles. When gates accumulate, infinitely many
can be crossed in a finite time and this count is useless. This section
uses a local exit estimate and one-way passage to prove Theorem~\ref{thm:A}.
We first give a quantitative version of the stopped-drift estimate for a
block with a gate at its entrance.

\subsection{Passage through a block}

\begin{proposition}[Uniform passage through a block]\label{prop:block}
Let $[h,d]\subset(l,r)$, $\sigma(h)=0$ and $f\ge\beta>0$ on $[h,d]$. No
assumption is made on the zeros of $\sigma$ in $(h,d)$. For $x\in[h,d)$,
\begin{equation}\label{eq:block-mean}
 T_d<\zeta\ \text{ almost surely},\qquad
 X_t\ge h\quad(0\le t\le T_d),\qquad
 \EE_xT_d\le\frac{d-x}{\beta}.
\end{equation}
If gates $h_n\uparrow c<d$ lie in $[h,d)$, then a path started at $h$ satisfies
$T_{h_n}\uparrow T_c<T_d<\infty$, and after $T_c$ it moves strictly to the
right of $c$. The reflected statement holds for $f\le-\beta<0$.
\end{proposition}

\begin{proof}
Since $f(h)>0$, Theorem~\ref{thm:one-way} keeps the path at or to the right of
$h$. Up to $T_d$ it therefore stays in the compact interval $[h,d]$, where it
cannot explode, and the stopped stochastic integral is a martingale. Hence
\[
 d\ge\EE_xX_{t\wedge T_d}
 =x+\EE_x\int_0^{t\wedge T_d}f(X_s)\,ds
 \ge x+\beta\,\EE_x(t\wedge T_d),
\]
and monotone convergence gives~\eqref{eq:block-mean}. By continuity every level
between $x$ and $d$ is reached before $T_d$, and the increasing times $T_{h_n}$
converge to $T_c$. The point $c$ is a zero of $\sigma$ with $f(c)\ge\beta$, hence
a rightward gate, and Theorem~\ref{thm:one-way} applies at $T_c$.
\end{proof}

The bound depends only on the width of the block and the drift bound $\beta$. At
an equilibrium, by contrast, accumulating zeros are passed at times that grow
without bound.

\begin{corollary}[Two kinds of accumulation points]\label{cor:cluster}
Every accumulation point $c$ of $Z$ inside $(l,r)$ is a gate. In a sufficiently
small neighborhood of $c$, Proposition~\ref{prop:block} applies to each block
whose entrance is a gate and whose orientation agrees with $f(c)$. In particular
one may start a block at $c$ itself. A block extending across $c$ requires a
gate on its incoming side. If zeros $h_n$ tend to a finite outer equilibrium
$e$ and a path reaches them in order, then $T_{h_n}\to\infty$.
\end{corollary}

\begin{proof}
By continuity $\sigma(c)=0$, and $f(c)\ne0$ because $(l,r)$ contains no
equilibrium; so $f$ has a fixed sign, bounded away from zero, near $c$. In the
equilibrium case, a finite limit of the times $T_{h_n}$ would give a finite
hitting time of $e$, contradicting Proposition~\ref{prop:inaccessible}.
\end{proof}

\begin{example}[Finite and infinite accumulation of passage times]\label{ex:clusters}
Let $\psi(z)=e^{-1/z^2}\sin(1/z)$ for $z\ne0$ and $\psi(0)=0$. This function is
smooth, bounded and globally Lipschitz, and its zeros $1/(n\pi)$ accumulate at
zero from both sides.

\emph{Gates accumulating at a gate.} Take $f\equiv1$ and $\sigma(x)=\psi(x-c)$.
There is no equilibrium, so $l=-\infty$ and $r=+\infty$. Every zero is a rightward
gate, and a path started at a zero $h<c$ passes all gates in $(h,c)$ before the
finite time $T_c$, with $\EE_hT_d\le d-h$ for each $d>c$. As $|z|\to\infty$,
$\psi(z)\sim1/z$, so $2f/\sigma^2\sim2(x-c)^2$. Hence $S(+\infty)<\infty$ on the
rightmost component and $S(-\infty)=-\infty$ on the leftmost. No component has
two repelling ends, and Theorem~\ref{thm:A} gives $X_t\to+\infty$ almost surely.
The coefficients are bounded, so there is no explosion.

\emph{Gates accumulating at an equilibrium.} Take $f(x)=-x$ and
$\sigma(x)=\psi(x)$. Now zero is an equilibrium and the zeros $h_n=1/(n\pi)$ are
leftward gates. From $x>1/\pi$, the positive solution is a nonnegative
supermartingale with $\EE_xX_t=xe^{-t}$, so $X_t\to0$ almost surely. Every gate
$h_n$ is reached, yet $T_{h_n}\to\infty$.
\end{example}

\subsection{The global alternative}

The proof of Theorem~\ref{thm:A} needs one more ingredient: an oscillating path
cannot settle at an interior point.

\begin{lemma}[Limits lie at the outer endpoints]\label{lem:outer-limit}
Almost surely, if $\lim_{t\uparrow\zeta}X_t$ exists in $[l,r]$, it equals $l$ or
$r$.
\end{lemma}

\begin{proof}
Each $z\in(l,r)$ has a compact neighbourhood $[p,q]\subset(l,r)$ on which either
$|\sigma|$ is bounded below by a positive constant or $f$ has a fixed sign and
$|f|\ge\beta>0$. In the first case $\sup_y\EE_y\tau_{(p,q)}<\infty$ by
Lemma~\ref{lem:regularity}; in the second,
$\beta\,\EE_y(t\wedge\tau_{(p,q)})\le q-p$ by the stopped drift identity. In
both cases the path leaves $(p,q)$ almost surely. Choose countably many such
neighbourhoods whose interiors cover $(l,r)$. A path with an interior limit is
bounded, so $\zeta=\infty$, and after some integer time it stays in the interior
of one of them. The Markov property at integer times and a countable union show
that this has probability zero.
\end{proof}

\begin{proof}[Proof of Theorem~\ref{thm:A}]
The path never reaches a finite outer equilibrium
(Proposition~\ref{prop:inaccessible}). Since the coefficients are locally
Lipschitz on a neighbourhood of each finite endpoint, explosion can only mean
$|X_t|\to\infty$, that is, convergence to an infinite outer endpoint.

\emph{Recurrent components are invariant.} Let $J\in\Rc$. An end of $J$ in $Z$
is repelling, so its drift points into $J$ (Lemma~\ref{lem:gate-ends}) and it is
never crossed from inside (Theorem~\ref{thm:one-way}). An equilibrium end is never
reached, and a repelling infinite end is not reached in finite time, since
\eqref{eq:feller} requires a finite scale limit. Theorem~\ref{thm:component}(i)
then gives recurrence, the $\liminf$ and $\limsup$, and the type of recurrence.

\emph{The alternative.} Let $L=\liminf_{t\uparrow\zeta}X_t$ and
$R=\limsup_{t\uparrow\zeta}X_t$ in $[l,r]$. If $L=R$, the path converges, and
Lemma~\ref{lem:outer-limit} gives alternative~(i). Suppose $L<R$. Then
$\zeta=\infty$, because explosion forces $|X_t|\to\infty$.

First, $(L,R)$ contains no point of $Z$. Indeed, if $z\in Z\cap(L,R)$, the path
is below and above $z$ at arbitrarily large times, so by continuity it visits
$z$ infinitely often, whereas a gate is visited at most once
(Theorem~\ref{thm:one-way}). This must hold simultaneously for all gates. The
set $Z$ has a countable dense subset $D$, and Theorem~\ref{thm:one-way} holds
simultaneously for all points of $D$. If $Z\cap(L,R)$ were nonempty, it would be
a nonempty relatively open subset of $Z$ and would contain a point of $D$. Hence
$(L,R)$ lies in a single component $J=(u,v)$.

Next, the path eventually stays in $J$. Pick a rational $c\in(L,R)$; it is
visited at arbitrarily large times. Leaving $J$ after a visit to $c$ would mean
crossing a gate end of $J$, after which $c$ is never visited again
(Theorem~\ref{thm:one-way}), or reaching an equilibrium or infinite end, which
is excluded. So $X_t\in J$ for all $t\ge T_c$.

Finally, $J\in\Rc$. Otherwise $J$ has an attracting end, and by
Theorem~\ref{thm:component}(ii)--(iii) and the strong Markov property at $T_c$,
a path that stays in $J$ after $T_c$ converges to an end of $J$ almost surely,
which contradicts $L<R$. There are countably many components and rationals, so
these null events can be discarded together. Thus alternative~(ii) holds with
the random time $T_c$, and the two alternatives are disjoint because $L=R$ in~(i)
and $L<R$ in~(ii).
\end{proof}

For finitely many gates, Theorem~\ref{thm:A} reduces to
Theorem~\ref{thm:routes}, and the route graph shows which component is selected.
For an arbitrary zero set there may be no graph to draw, since a block can hold
infinitely many vertices. One-way passage excludes a gate between the lower
and upper limits of an oscillating path. Local exit estimates exclude an
interior accumulation point as a path limit. An outer equilibrium can still
be a limit, approached only as $t\to\infty$.

\section{Stationary and limiting laws}\label{sec:laws}

Theorem~\ref{thm:A} splits the path space into the events
\[
 E_l=\{X_t\to l\},\qquad E_r=\{X_t\to r\},\qquad
 E_J=\{X_t\in J\ \text{for all large } t\}\quad(J\in\Rc),
\]
with probabilities $q_l(x)$, $q_r(x)$ and $w_J(x)$ summing to one. This section
turns that pathwise picture into statements about laws. Positive recurrent
components contribute their stationary densities. A null recurrent component
with infinite speed mass at just one end contributes a point mass there,
even though its paths keep returning. After an explosion we set $X_t$ equal to the infinite endpoint
reached, for $t\ge\zeta$, and we work on the two-point compactification of
$[l,r]$.

\subsection{Stationary laws}

\begin{proposition}[All stationary laws]\label{prop:stationary}
Theorem~\ref{thm:B}\textup{(i)} holds.
\end{proposition}

\begin{proof}
The point mass at a finite equilibrium is stationary, and so is $\pi_J$ for
$J\in\Rc_+$, because $J$ is invariant (Theorem~\ref{thm:A}) and
Lemma~\ref{lem:recurrence} applies. Conversely, let $\nu$ be a stationary
probability measure, and let $D=(l,r)\setminus\bigcup_{J\in\Rc}J$. For each
compact $K\subset(l,r)$, every path eventually avoids $K\cap D$: on $E_l\cup E_r$
it leaves $K$, and on $E_J$ it stays in $J$. Bounded convergence gives
$\frac1t\int_0^t\PP_y(X_s\in K\cap D)\,ds\to0$ for each $y$. Stationarity and
dominated convergence then give
\[
 \nu(K\cap D)=\int\frac1t\int_0^t\PP_y(X_s\in K\cap D)\,ds\,\nu(dy)\longrightarrow0,
\]
so $\nu(D)=0$. Each $J\in\Rc$ is invariant, so the restriction of $\nu$ to $J$
is stationary for the process on $J$. By Lemma~\ref{lem:recurrence} it is zero
when $J$ is null recurrent and a multiple of $\pi_J$ when $J\in\Rc_+$. The
remaining mass sits at the finite outer equilibria.
\end{proof}

\subsection{Null recurrent components}

In a null recurrent component the path visits every level infinitely often, yet
it spends a vanishing fraction of time in any set of finite speed mass.

\begin{lemma}[Escape from finite-speed sets]\label{lem:null-escape}
Let $J\in\Rc$ have $M_J=\infty$, and let $A\subset J$ be a Borel set with
$\int_Am_J<\infty$. Then, for every $x\in J$,
\[
 \PP_x(X_t\in A)\longrightarrow0,\qquad
 \frac1t\int_0^t\mathbf1_A(X_s)\,ds\longrightarrow0\quad\PP_x\text{-a.s.}
\]
\end{lemma}

The proof uses the symmetric transition density with respect to the speed
measure, together with a cycle decomposition; it is given in
Appendix~\ref{app:null}. The interval $J$ may be unbounded.
For uniformly nondegenerate coefficients on $\R$, the first statement is the
sweeping property of~\cite[Lemma~1]{Rudnicki1993}, which holds exactly when the
speed mass is infinite.

\begin{theorem}[Limits in a null recurrent component]\label{thm:null}
Let $J=(u,v)\in\Rc$ be null recurrent with infinite speed mass near $u$ only,
that is, $\int_u^cm_J=\infty$ and $\int_c^vm_J<\infty$ for $c\in J$. Then, on the
two-point compactification of $[u,v]$,
\[
 \mathcal L_x(X_t)\Longrightarrow\delta_u,\qquad
 \frac1t\int_0^t\delta_{X_s}\,ds\Longrightarrow\delta_u\quad\PP_x\text{-a.s.},
\]
while $\liminf X_t=u$ and $\limsup X_t=v$ almost surely. For $u=-\infty$,
$\PP_x(X_t\ge c)\to0$ for each $c\in J$. The conclusion for a unique
infinite-speed right end follows by reflection. If both ends have infinite
speed mass, every subsequential weak limit of the laws, and almost
surely of the empirical measures, is carried by $\{u,v\}$.
\end{theorem}

\begin{proof}
For each $c\in J$, the set $[c,v)$ has finite speed mass, so by
Lemma~\ref{lem:null-escape} its probability and occupation fraction tend to zero.
Taking a sequence $c\downarrow u$ and using continuity of test functions on the
compactified interval proves both limits. The $\liminf$ and $\limsup$ come from
recurrence. In the two-sided case, apply Lemma~\ref{lem:null-escape} to an
exhaustion of $J$ by compact intervals.
\end{proof}

The two-sided statement cannot be improved in general: for
$dX=dW+b(X)\,dt$ with $b(x)\to0$ as $|x|\to\infty$, Rudnicki~\cite{Rudnicki1993}
constructs $b$ so that both $\delta_{-\infty}$ and $\delta_{+\infty}$ occur as
subsequential limits of the laws.

\begin{example}[A law tending to an unstable equilibrium]\label{ex:null-dirac}
On $(0,1)$ let $\sigma(x)=x(1-x)$ and $f(x)=\tfrac12x(1-x)(1-3x)$. Both ends are
equilibria, and one may take $\rho(x)=x^{-1}(1-x)^{-2}$ and $m(x)=2/x$. Both
scale limits are infinite, but the speed mass diverges only at zero. The process
is null recurrent, $\mathcal L_x(X_t)\Rightarrow\delta_0$, and the empirical
measures converge to $\delta_0$ almost surely. Nevertheless $\liminf X_t=0$ and
$\limsup X_t=1$. Zero is unstable in probability: convergence of the law to it
does not give the uniform-in-time control of Definition~\ref{def:stable}.
\end{example}

\begin{corollary}[Where null recurrence can occur]\label{cor:null-location}
Suppose $Z\ne\varnothing$. A null recurrent component touches exactly one outer
endpoint $e_J\in\{l,r\}$, and its speed mass is infinite near $e_J$ and finite
near its other end. There are at most two null recurrent components.
\end{corollary}

\begin{proof}
An end of $J\in\Rc$ that lies in $Z$ is a gate with inward drift, near which the
speed mass is finite (Lemma~\ref{lem:gate-ends}). The speed mass is also finite
on compact subintervals. Since $Z\ne\varnothing$, no component equals $(l,r)$, so
each component has at most one outer end. A null recurrent component therefore
has exactly one, carrying all of the infinite mass, and each outer endpoint
belongs to at most one component.
\end{proof}

\subsection{Proof of Theorem~\ref{thm:B}}

Part~(i) is Proposition~\ref{prop:stationary}. For~(ii) and~(iii), fix
$x\in(l,r)$ and assume $Z\ne\varnothing$. Let $\Theta$
be the first integer $n$ with $X_n\in\bigcup_{J\in\Rc}J$. On $E_J$, $\Theta$ is
finite and $X_\Theta\in J$, and the path stays in $J$ afterwards.

\emph{Laws.} Let $g$ be continuous on the compactified interval. On $E_l$ and
$E_r$, $g(X_t)$ converges to $g(l)$ or $g(r)$ along every path. On $E_J$,
conditioning on $\Theta=n$ and $X_n=y$, the Markov property gives
$\EE_y\,g(X_{t-n})$. This tends to $\int g\,d\pi_J$ if $J\in\Rc_+$
(Lemma~\ref{lem:recurrence}) and to $g(e_J)$ if $J$ is null recurrent
(Theorem~\ref{thm:null} and Corollary~\ref{cor:null-location}). Dominated
convergence over $n$, $y$ and the countably many $J$ gives~\eqref{eq:weak-limit}.

\emph{Total variation.} At time $t$, the law gives mass $\PP_x(\zeta\le t)$ to
the outer endpoints. The limit gives them mass $Q_l(x)+Q_r(x)$. Hence the
total-variation distance is at least $Q_l(x)+Q_r(x)-\PP_x(\zeta<\infty)$. If this
difference vanishes, the event of selecting a null recurrent component and
the event of approaching an endpoint without explosion both have probability
zero. On the explosion event couple $X_t$ to the endpoint actually reached:
their discrepancy has probability at most $\PP_x(t<\zeta<\infty)$.
On $\{\Theta=n,X_n\in J\}$ with $J\in\Rc_+$, for $t\ge n$ the conditional
kernel from $X_n$ converges to $\pi_J$ in total variation by
Lemma~\ref{lem:recurrence}. Integrate its bounded total-variation distance
over $X_n$ and sum over $n,J$, bounding the contribution with $n>t$ by
$\PP_x(t<\Theta<\infty)$. This proves sufficiency.

\emph{Empirical measures.} Along paths in $E_l$ or $E_r$ the empirical measures
converge to the point mass at the limit. On $E_J$ with $J\in\Rc_+$, the ergodic
theorem of Lemma~\ref{lem:recurrence}, applied after $\Theta$, gives $\pi_J$; for
null recurrent $J$, Theorem~\ref{thm:null} gives $\delta_{e_J}$. Using a countable
dense family of test functions gives a single event of probability one, and the
initial segment before $\Theta$ does not affect time averages.
\qed

\section{Quantitative passage and selection}\label{sec:times}

On a finite route, the selection time $\tau_*$ is the time of the last gate
crossing, with $\tau_*=0$ if there is none. Theorem~\ref{thm:routes} gives
almost-sure finiteness. We estimate its exponential moment for gate-bounded
stages. We also compute the mean passage time from a fixed point between a
repelling equilibrium and the first outward gate. For an arbitrary zero set,
the block estimate below applies whether or not a last crossing exists.

\begin{lemma}[Geometric tails]\label{lem:geometric}
Let $\tau$ be the first hitting time of a fixed target, and suppose the path stays
in a set $D$ before $\tau$, with $\sup_{y\in D}\EE_y\tau\le C$ for some
$0<C<\infty$. Then, for each integer $k\ge0$,
\[
 \sup_{y\in D}\PP_y(\tau>2kC)\le2^{-k},\qquad
 \sup_{y\in D}\EE_ye^{\theta\tau}\le\frac{e^{2C\theta}}{2-e^{2C\theta}}
 \quad\Big(0<\theta<\frac{\log2}{2C}\Big).
\]
\end{lemma}

\begin{proof}
Markov's inequality gives $\PP_y(\tau>2C)\le\tfrac12$ uniformly in $y\in D$. The
Markov property at the times $2C,4C,\dots$ gives the tail bound. So
$\lceil\tau/2C\rceil$ is dominated by a geometric variable with parameter
$\tfrac12$, and summing its exponential moments gives the second bound.
\end{proof}

\begin{proof}[Proof of Theorem~\ref{thm:C}\textup{(i)}]
Proposition~\ref{prop:block} gives $\EE_yT_d\le(d-h)/\beta$ for all $y\in[h,d)$,
and the path stays in $[h,d)$ before $T_d$. Apply Lemma~\ref{lem:geometric} with
$C=(d-h)/\beta$.
\end{proof}

\begin{proof}[Proof of Theorem~\ref{thm:C}\textup{(ii)}]
On leaving its initial component through a gate, the path enters the next
component at an inward gate. At each gate-bounded passage stage take the
target to be the set of its outward gates, allowing two possible exits at
the initial stage. Proposition~\ref{prop:mean-exit} bounds the mean time to
this target uniformly over the interior and any inward entrance gate.
Lemma~\ref{lem:geometric} thus bounds the conditional exponential moment
of the stage duration by $B_C(\theta)=e^{2C\theta}/(2-e^{2C\theta})$.
Iterate at the successive gate hitting times by the strong Markov property.
A completed route contributes zero duration in every unused stage. There
are at most $|Z|$ stages of positive duration, giving the stated product.
\end{proof}

The same product argument applies to a prescribed finite sequence of passage
blocks, using a factor for each stage with a uniform mean bound. An approach
to a repelling equilibrium requires a separate calculation because the initial
point in the adjacent component can be arbitrarily close to that equilibrium.

\begin{proposition}[The stage next to an equilibrium]\label{prop:eq-stage}
Assume the setting of Theorem~\ref{thm:C}\textup{(iii)}. Then $T_h<\infty$
almost surely, and
\begin{equation}\label{eq:eq-stage}
 \EE_xT_h=\int_l^h\big(B-S(x\vee y)\big)\,m(y)\,dy
 =(B-S(x))\int_l^xm+\int_x^h(B-S(y))\,m(y)\,dy,
\end{equation}
with the value $+\infty$ allowed. In particular $\EE_xT_h<\infty$ if and only if
$\int_lm<\infty$ near $l$. If $Z$ is finite and every later stage is bounded by
two gates, then $\EE_x\tau_*<\infty$ under the same condition.
\end{proposition}

\begin{proof}
Theorem~\ref{thm:component}(ii) gives $T_h<\infty$ almost surely. For $l<d<x$, the
interval $(d,h)$ has a regular end $d$ and an outward gate end $h$.
Lemma~\ref{lem:green}, obtained by exhaustion with regular barriers, gives
\[
 \EE_x(T_d\wedge T_h)=\int_d^hG_d(x,y)\,m(y)\,dy,\qquad
 G_d(x,y)=\frac{\big(S(x\wedge y)-S(d)\big)\big(B-S(x\vee y)\big)}{B-S(d)} .
\]
As $d\downarrow l$, $S(d)\to-\infty$, so $G_d(x,y)$ increases to $B-S(x\vee y)$.
Also $T_d\wedge T_h$ increases to $T_h$, because $l$ is never reached. Monotone
convergence gives~\eqref{eq:eq-stage}. The second integral in~\eqref{eq:eq-stage}
is finite by Lemma~\ref{lem:gate-ends}(ii), and $B-S(x)>0$, which gives the
criterion. The last claim follows from Theorem~\ref{thm:C}(ii) and the strong
Markov property at $T_h$.
\end{proof}

\begin{proof}[Proof of Theorem~\ref{thm:C}\textup{(iii)--(iv)}]
Part~(iii) is Proposition~\ref{prop:eq-stage}. For~(iv), let $V(x)=\log(x-l)$.
Then
\[
 \cL V(x)=\frac{f(x)}{x-l}-\frac{\sigma^2(x)}{2(x-l)^2}\longrightarrow\lambda
 \qquad(x\downarrow l),
\]
so $\cL V\ge c_0$ on $(l,d)$ for $d$ close to $l$. Stopping at
$T_d\wedge T_{l+\varepsilon}\wedge t$ and using $V\le\log(d-l)$ there, Itô's
formula gives
$c_0\,\EE_x(t\wedge T_d\wedge T_{l+\varepsilon})\le\log\big((d-l)/(x-l)\big)$. Let
$\varepsilon\downarrow0$, using that $l$ is never reached, and then let
$t\to\infty$.
\end{proof}

Strict repulsion also gives finite speed mass near $l$, including when
$\sigma'(l)=0$. Indeed, for $l<x<d$ as in~(iv), Lemma~\ref{lem:green} on
$(l,d)$ gives $\EE_xT_d\ge (S(d)-S(x))\int_l^x m(y)\,dy$. The coefficient
is positive, and~(iv) bounds the left side for each fixed $x$.

\begin{example}[Almost-sure selection with infinite mean]\label{ex:infinite-mean}
On $[0,2]$ take the locally Lipschitz coefficients
\[
 \sigma(x)=\begin{cases}x(1-x),&0\le x\le1,\\(x-1)(2-x),&1\le x\le2,\end{cases}
 \qquad
 f(x)=\begin{cases}x/2,&0\le x\le1,\\(2-x)(3-\tfrac52x),&1\le x\le2,\end{cases}
\]
extended to a neighbourhood of $[0,2]$. The points $0$ and $2$ are consecutive
equilibria, and $1$ is a rightward gate with $f(1)=\tfrac12$. At zero,
$\lambda=0$ with smooth coefficients, so zero is repelling with infinite speed
mass (Corollary~\ref{cor:first-order}): $\rho(x)\sim C/x$ and $m(x)\sim D/x$. At
two, $f'(2)=2$ and $\sigma'(2)=-1$, so $\lambda=\tfrac32>0$: two is repelling with
finite speed mass. The component $(1,2)$ is terminal and positive recurrent. For
$x\in(0,1)$, $\tau_*=T_1<\infty$ almost surely, but $\EE_x\tau_*=\infty$ by
Proposition~\ref{prop:eq-stage}. Neither almost-sure selection nor positive
recurrence of the selected component implies a finite mean selection time.
\end{example}

\section{Discussion}\label{sec:discussion}

Two mechanisms govern a scalar SDE whose noise vanishes. Equilibria are never
reached from another state, and a gate can be crossed at most once.
For an arbitrary closed set $E=\{f=\sigma=0\}$, a start in $E$ is stationary.
Every other initial state belongs to a connected component of $\R\setminus E$.
Each finite endpoint of that component belongs to $E$, so
Theorem~\ref{thm:A} applies there even when equilibria accumulate.

\begin{enumerate}
\item \emph{Two infinite-speed ends.} In the null recurrent case with
$Z=\varnothing$ and infinite speed mass at both ends,
Theorem~\ref{thm:null} places all subsequential limits on the two endpoints.
For uniformly nondegenerate, bounded coefficients on $\R$, Rudnicki~\cite{Rudnicki1993}
gives conditions under which the weight at $+\infty$ converges, with an explicit
limit, and constructs an equation for which even its time averages oscillate
between $0$ and $1$. Extending such criteria to components whose ends are
degenerate equilibria remains open.
\item \emph{Rates.} We give no rate of convergence to $\pi_J$, and no uniform
selection-time bounds along infinite routes without a uniform drift bound.
\item \emph{Higher dimensions.} A natural analogue of a gate is a hypersurface on
which the noise degenerates while the drift crosses it transversally. Whether a
similar global alternative holds there is open.
\end{enumerate}

\appendix
\section{Classical tools}\label{app:tools}

This appendix collects the classical facts used above, with short proofs or
references.

\subsection{Regularity on a nondegenerate interval}\label{app:regularity}

\begin{lemma}\label{lem:regularity}
Let $[a,b]\subset(l,r)$ with $|\sigma|\ge\kappa>0$ and $|f|\le F$ on $[a,b]$.
Then $\sup_{x\in[a,b]}\EE_x\tau_{(a,b)}<\infty$, the two-barrier
formula~\eqref{eq:two-barrier} holds, and every level of a component is reached
with positive probability from every other level of that component, before the
path leaves it.
\end{lemma}

\begin{proof}
Choose $k>2F/\kappa^2$ and put $V(y)=(e^{kb}-e^{ky})/C$ with
$C=e^{ka}k(k\kappa^2/2-F)>0$. Then $V\ge0$ on $[a,b]$ and
$\cL V=-e^{ky}(kf+\tfrac12k^2\sigma^2)/C\le-1$, so stopped Itô calculus gives
$\EE_x\tau_{(a,b)}\le V(x)$. Since $\cL S=0$, the bounded stopped martingale
$S(X_{t\wedge\tau_{(a,b)}})$ gives~\eqref{eq:two-barrier}, and both exit
probabilities are positive.
\end{proof}

\subsection{Lyapunov criteria}\label{app:lyapunov}

We use the inward coordinate~\eqref{eq:inward} at an equilibrium, with generator
$\cL_s=\mu\partial_y+\tfrac12\nu^2\partial_{yy}$. The following criteria are
classical~\cite{Khasminskii2012}. Each is proved by Itô's formula stopped on
$[a,b]\subset(0,h)$, followed by $a\downarrow0$.

\begin{proposition}\label{prop:lyapunov}
\begin{enumerate}[label=\textup{(\roman*)}]
\item If $v\in C([0,h])\cap C^2((0,h))$, $v(0)=0$, $v>0$ on $(0,h]$ and
$\cL_sv\le0$, then $\PP_{e+sy}(H_b<\infty)\le v(y)/v(b)$ for $0<y<b<h$, so $e$ is
stable from side $s$.
\item If $v\in C^2((0,h))$ is positive, $v(y)\to\infty$ as $y\downarrow0$, and
$\cL_sv\le0$, then $e$ is unstable from side $s$.
\item With $R(x)=2(x-e)f(x)/\sigma^2(x)$: if $\limsup_{x\to e}R<1$ on side $s$,
then $e$ is stable from that side, using $v=y^p$ with small $p>0$; if
$\liminf_{x\to e}R>1$, it is unstable, using $v=y^{-p}$.
\end{enumerate}
\end{proposition}

Part~(i) with $v(y)=y^p$ and $0<p<\tfrac12$ gives Example~\ref{ex:noise-stability}.

\subsection{Green functions and mean exit times}\label{app:green}

\begin{lemma}\label{lem:green}
Let $J=(u,v)$ be an interval on which $\sigma\ne0$, assume at least one
end has finite scale, and let $\tau_J$ be its exit time, including the
lifetime if explosion occurs before any finite exit. With $G_J$ as in
Proposition~\ref{prop:mean-exit}, for every $x\in J$,
\[
 \EE_x\tau_J=\int_u^v G_J(x,y)m(y)\,dy
\]
as an identity in $[0,\infty]$.
If both scale limits are finite, write $A=S(u+)$ and $B=S(v-)$.
The mean is uniformly bounded in $x$ if $S(y)-A$ is integrable against $m$
near $u$ and $B-S(y)$ is integrable against $m$ near $v$.
If $S(u+)=-\infty$ and $B=S(v-)<\infty$, then
\[
 \sup_{x\in J}\EE_x\tau_J=\int_u^v (B-S(y))\,m(y)\,dy .
\]
For $S(u+)=A>-\infty$ and $S(v-)=+\infty$, the reflected supremum is
$\int_u^v(S(y)-A)m(y)\,dy$.
\end{lemma}

\begin{proof}
First consider regular barriers $u<a<x<b<v$. Write
$F(z)=S(z)-S(a)$, $R(z)=S(b)-S(z)$ and $D=S(b)-S(a)$. The kernel on $(a,b)$ is
\[
 G_{a,b}(x,y)=\frac{F(x\wedge y)R(x\vee y)}{D}.
\]
Its integral $U(x)=\int_a^bG_{a,b}(x,y)m(y)\,dy$ vanishes at $a,b$.
Splitting at $x$ and differentiating gives $(U'/\rho)'=-m$, hence
$\cL U=-1$. The stopped It\^o formula and the finite exit time of
Lemma~\ref{lem:regularity} give $U(x)=\EE_x(T_a\wedge T_b)$.

Take $a\downarrow u$, $b\uparrow v$ along nested intervals. The exit times
increase to $\tau_J$ by continuity and local continuation. The kernels,
extended by zero outside $(a,b)$, increase as either barrier moves outward:
for fixed $x,y$ this follows directly from the displayed fraction. Their
pointwise limit is the corresponding kernel in Proposition~\ref{prop:mean-exit}.
Monotone convergence gives the identity, including infinite means.

When both scale limits are finite, $G_J(x,y)\le
\min\{S(y)-A,B-S(y)\}$ for every $x,y$, giving the stated uniform bound.
With $u$ repelling and $v$ attracting, $G_J(x,y)=B-S(x\vee y)$ increases to
$B-S(y)$ as $x\downarrow u$. Monotone convergence gives the exact supremum.
Reflection gives the left-attracting case.
\end{proof}

\subsection{Recurrence and ergodicity}\label{app:recurrence}

Lemma~\ref{lem:recurrence} is classical~\cite[Ch.~15]{Karlin}. Recurrence follows
from~\eqref{eq:two-barrier} by letting one barrier tend to a repelling end, and
from the strong Markov property. On compact subintervals the killed transition
kernel has a strictly positive continuous density. Hence a fixed-time skeleton
chain is irreducible and aperiodic, and compact sets are small. When
$M_J<\infty$, the process is positive Harris recurrent, and total-variation
convergence and the ergodic theorem follow from~\cite{MeynTweedie1993}.

\subsection{Proof of Lemma~\ref{lem:null-escape}}\label{app:null}

\emph{Laws.} The transition semigroup of the diffusion on $J$ is self-adjoint on
$L^2(\mu)$, $\mu(dy)=m_J(y)\,dy$, with a continuous symmetric density
$p_t(x,y)$~\cite[Ch.~V]{RogersWilliams}. Its Dirichlet form is
$\mathcal E(g,g)=\int_J|g'|^2/\rho_J$, whose only zero-energy functions are
constants. Since $\mu(J)=\infty$, no nonzero constant lies in $L^2(\mu)$, and the
spectral theorem gives $\|P_tg\|_{L^2(\mu)}\to0$ for every $g\in L^2(\mu)$. For
$s>0$ fixed and $t>s$,
\[
 P_t\mathbf1_A(x)\le\|p_s(x,\cdot)\|_{L^2(\mu)}\,\|P_{t-s}\mathbf1_A\|_{L^2(\mu)}
 =\sqrt{p_{2s}(x,x)}\,\|P_{t-s}\mathbf1_A\|_{L^2(\mu)}\longrightarrow0 .
\]
\emph{Occupation times.} From any initial $x\in J$ the hitting time of
$c\in J$ is almost surely finite, and this initial segment has vanishing
weight in the time average. Fix $d>c$ in $J$. Starting from $c$, a cycle
hits $d$ and returns to $c$. By recurrence, each cycle is finite. Its expected occupation
density is the constant multiple $(S_J(d)-S_J(c))\,m_J$, obtained by adding the
Green kernels of the two halves of the cycle. So cycles have infinite mean
length, but the time spent in $A$ during a cycle has finite mean
$(S_J(d)-S_J(c))\mu(A)$. Cycles are independent and identically distributed by
the strong Markov property. The strong law, applied to truncated cycle lengths,
gives $N(t)/t\to0$ for the number $N(t)$ of cycles completed by time $t$. The
strong law for the occupation rewards bounds the time spent in $A$ by the
sum of the first $N(t)+1$ full rewards, including the cycle in progress.
This sum is $O(N(t)+1)=o(t)$ almost surely.
\section*{Disclosure of AI usage}

During manuscript writing, the authors used OpenAI's Codex to review theorem
statements and proofs, identify missing hypotheses, examine counterexamples,
and suggest revisions to the mathematical arguments and exposition.
The authors independently verified the mathematical arguments in the final manuscript and take responsibility for its contents.

\end{document}